\documentclass[12pt]{article}

\usepackage{latexsym,amsfonts,amsmath,amssymb,graphicx,hyperref}

\newcommand{\comment}[1]{}

\newcommand{\St}{\mathbb{S}}
\newcommand{\EE}{\mbox{\bf E}\,}
\newcommand{\PP}{\mbox{\bf P}\,}
\newcommand{\QQ}{\mbox{\bf Q}\,}

\newcommand{\R}{\mathbb{R}}
\newcommand{\C}{\mathbb{C}}
\newcommand{\Q}{\mathbb{Q}}
\newcommand{\HH}{\mathbb{H}}
\newcommand{\N}{\mathbb{N}}

\newcommand{\B}{\mbox{\bf B}}

\newcommand{\pa}{\partial}

\newcommand{\F}{{\cal F}}

\newcommand{\no}{\noindent}

\newcommand{\BGE}{\begin{equation}}
\newcommand{\BGEN}{\begin{equation*}}
\newcommand{\EDE}{\end{equation}}
\newcommand{\EDEN}{\end{equation*}}

\def\eps{\varepsilon}
\def\til{\widetilde}
\def\ha{\widehat}
\def\sem{\setminus}

\def\lin{\overline}
\def\vphi{\varphi}

\DeclareMathOperator{\sign}{sign} 
\DeclareMathOperator{\dist}{dist} 
\DeclareMathOperator{\hcap}{hcap} 
\DeclareMathOperator{\Imm}{Im } \DeclareMathOperator{\Ree}{Re }

\DeclareMathOperator{\scap}{scap} 
\DeclareMathOperator{\out}{out}

\def\h0{{\bf h}}

\newtheorem{Lemma}{Lemma}[section]
\newtheorem{Theorem}{Theorem}[section]

\newtheorem{Corollary}{Corollary}[section]
\newtheorem{Proposition}{Proposition}[section]
\newtheorem{Conjecture}{Conjecture}
\numberwithin{equation}{section}

\begin{document}
\title{\bf Duality of Chordal SLE, II}
\date{\today}
\author{Dapeng Zhan\footnote{Yale University}}
\maketitle
\begin{abstract}
We improve the geometric properties of  SLE$(\kappa;\vec{\rho})$
processes derived in an earlier paper, which are then used to obtain
more results about the duality of SLE. We find that for $\kappa\in (4,8)$, the
boundary of a standard chordal SLE$(\kappa)$ hull stopped on swallowing
a fixed $x\in\R\sem\{0\}$ is the image of some SLE$(16/\kappa;\vec{\rho})$ trace
started from a random point. Using this fact together with a similar proposition
in the case that $\kappa\ge 8$, we obtain a description of the boundary of
a standard chordal SLE$(\kappa)$ hull for $\kappa>4$, at a finite stopping time.
Finally, we prove that for $\kappa>4$,  in many cases,  the limit of
a chordal or strip SLE$(\kappa;\vec{\rho})$ trace exists.
\end{abstract}

\section{Introduction}
This paper is a follow-up of the paper \cite{duality}, in which we proved some versions of
Duplantier's duality conjecture about Schramm's SLE process (\cite{S-SLE}). In the present paper, we will improve the technique
used in \cite{duality}, and derive more results about the duality conjecture.

Let us now briefly review some results in \cite{duality}. Let $\kappa_1<4<\kappa_2$ with $\kappa_1\kappa_2=16$.
Let $x_1\ne x_2\in\R$. Let $N\in\N$ and $p_1,\dots,p_N\in\R\cup\{\infty\}\sem\{x_1,x_2\}$ be distinct points. Let $C_1,\dots,C_N\in\R$ and $\rho_{j,m}=C_m(\kappa_j-4)$,
$1\le m\le N$, $j=1,2$. Let $\vec{p}=(p_1,\dots,p_N)$ and $\vec{\rho}_j=(\rho_{j,1},\dots,\rho_{j,N})$, $j=1,2$.
Using the method of coupling two SLE processes obtained
in \cite{reversibility} and some computations about SLE$(\kappa;\vec{\rho})$ processes, we derived
Theorem 4.1 in \cite{duality}, which says that there is a coupling of a chordal SLE$(\kappa_1;-\frac{\kappa_1}2,\vec{\rho}_1)$
process $K_1(t)$, $0\le t<T_1$, started from $(x_1;x_2,\vec{p})$, and a chordal SLE$(\kappa_2;-\frac{\kappa_2}2,\vec{\rho}_2)$
process $K_2(t)$, $0\le t<T_2$, started from $(x_2;x_1,\vec{p})$, such that certain properties are satisfied.
Moreover, some $p_m$ could take value $x_j^\pm$, $j=1,2$, if the corresponding force $\rho_{j,m}\ge\kappa_j/2-2$.

This theorem was then applied to the case that $N=3$; $x_1<x_2$; $p_1\in(-\infty,x_1)$
or $=x_1^-$; $p_2\in(x_2,\infty)$, or $=\infty$, or $=x_2^+$; and $p_3\in(x_1,x_2)$, or $=x_1^+$,
or $=x_2^-$; $C_1\le 1/2$, $C_2=1-C_1$, and $C_3=1/2$. Using some geometric properties about SLE$(\kappa;\vec{\rho})$ processes, we concluded
that $K_1(T_1^-):=\cup_{0\le t<T_1} K_1(t)$  is
the outer boundary of $K_2(T_2^-):=\cup_{0\le t<T_2} K_2(t)$ in $\HH$.

The following proposition, i.e., Theorem 5.2 in \cite{duality}, is an application of the above result. It describes the
boundary of a standard chordal SLE$(\kappa)$ hull, where $\kappa\ge 8$, at the time when a fixed $x\in\R\sem\{0\}$ is swallowed.

\begin{Proposition} Suppose $\kappa\ge 8$, and $K(t)$, $0\le t<\infty$,
is a standard chordal SLE$(\kappa)$ process. Let $x\in\R\sem\{0\}$ and
$T_x$ be the first $t$ such that $x\in\lin{K(t)}$. Then $\pa
K(T_x)\cap\HH$ has the same distribution as the image of a chordal
SLE$(\kappa';-\frac{\kappa'}2,-\frac{\kappa'}2,\frac{\kappa'}2-2)$
trace started from $(x;0,x^a,x^b)$, where $\kappa'=16/\kappa$, $a=\sign(x)$
and $b=\sign(-x)$. So a.s.\ $\pa K(T_x)\cap\HH$ is a crosscut in $\HH$ on $\R$ connecting $x$
with some $y\in\R\sem\{0\}$ with $\sign(y)=\sign(-x)$.
\label{kappa>8}
\end{Proposition}

Here a crosscut in $\HH$ on $\R$ is a simple curve in $\HH$ whose two ends approach to two different points on $\R$.
Since $\kappa\ge 8$, the trace is space-filling, so a.s.\ $x$ is visited by the trace at time $T_x$, and so $x$ is an end point of
$\lin{K(T_x)}\cap\R$. From this proposition, we see that the boundary of $K(T_x)$ in $\HH$ is
an SLE$(16/\kappa)$-type trace in $\HH$ started from $x$.

\vskip 3mm

The motivation of the present paper is to derive the counterpart of Proposition \ref{kappa>8} in the
case that $\kappa\in(4,8)$. In this case, the trace, say $\gamma$, is not space-filling, so
a.s.\ $x$ is not visited by $\gamma$, at time $T_x$, and so $x$ is an interior point of $\lin{K(T_x)}\cap\R$.
Thus we can not expect that the boundary of $K(T_x)$ in $\HH$ is a curve started from $x$.

This difficulty will be overcome by conditioning the process $K(t)$, $0\le t<T_x$, on the value of $\gamma(T_x)$.
The conditioning should be done carefully since the probability that $\gamma(T_x)$ equals to any particular value is zero. 
Instead of taking limits, we will express $K(t)$, $0\le t<T_x$, as an integration
of some SLE$(\kappa;\vec{\rho})$ processes. 
In Section \ref{SLE-intergration}, we will prove that the distribution of $K(t)$, $0\le t<T_x$, is an integration of 
the distributions of SLE$(\kappa;-4,\kappa-4)$ processes started from $(0;y,x)$ against $d\lambda(y)$, where $\lambda$
is the distribution of $\gamma(T_x)$. This is the statement of Corollary \ref{kappa(4,8)}.

In Section \ref{GeometricProperties}, we will improve the geometric results about SLE$(\kappa;\vec{\rho})$ processes
that were derived in \cite{duality}. Using these geometric results, we will prove in Section \ref{Duality} that Proposition \ref{bartk}
can be applied with $N=4$ and suitable values of $p_m$ and $C_m$ for $1\le m\le 4$, to obtain more results about duality.
Especially, using Corollary \ref{kappa(4,8)}, we will obtain the
counterpart of Proposition \ref{kappa>8} in the case that $\kappa\in(4,8)$, which is Theorem \ref{kappa(4,8)swallow} below.

\begin{Theorem} Let $\kappa\in(4,8)$, and $x\in\R\sem\{0\}$. Let $K(t)$ and $\gamma(t)$,
$0\le t<\infty$, be  standard chordal SLE$(\kappa)$ process and trace, respectively.
Let $T_x$ be the first time that $x\in \lin{K(t)}$. Let $\bar\mu$
denote the distribution of
$\pa K(T_x)\cap\HH$. Let $\lambda$ denote the distribution of $\gamma(T_x)$. Let $\kappa'=16/\kappa$,
$a=\sign(x)$ and $b=\sign(-x)$. Let $\bar\nu_y$ denote the distribution of the image of a
chordal SLE$(\kappa'; -\frac{\kappa'}2,\frac{3}2 \kappa'-4,-\frac{\kappa'}2+2,\kappa'-4)$ trace
started from $(y;0,y^a,y^b,x)$. Then  $\bar\mu=\int \bar\nu_y\,d\lambda(y)$. So a.s.\
$\pa K(T_x)\cap\HH$ is a crosscut in $\HH$ on $\R$ connecting some $y,z\in\R\sem\{0\}$, where $\sign(y)=\sign(x)$,
$|y|>|x|$, and $\sign(z)=\sign(-x)$.
\label{kappa(4,8)swallow}
\end{Theorem}

In Section \ref{BoundaryofChordalSLE}, we will use Theorem \ref{kappa(4,8)swallow} and Proposition \ref{kappa>8}
to study the boundary of a standard chordal SLE$(\kappa)$ hull, say $K(t)$, at a finite positive stopping time $T$.
Let $\gamma(t)$ be the corresponding SLE trace. We will find that if $\gamma(T)\in\R$, then
$\pa K(T)\cap\HH$ is a crosscut in $\HH$ with $\gamma(T)$ as one end point; and if $\gamma(T)\in\HH$, then $\pa K(T)\cap\HH$ is
the union of two semi-crosscuts in $\HH$, which both have $\gamma(T)$ as one end point. Here a semi-crosscut in $\HH$ is a
simple curve in $\HH$ whose one end lies in $\HH$ and the other end approaches to a point on $\R$.
Moreover, in the latter case, every intersection point
of the two semi-crosscuts other than $\gamma(T)$ corresponds to a cut-point of $K(T)$. If $\kappa\ge 8$, then the two
semi-crosscuts only meet at $\gamma(T)$, and so $\pa K(T)\cap\HH$ is again a crosscut in $\HH$ on $\R$.

In the last section of this paper, we will use the results in Section \ref{BoundaryofChordalSLE} to derive more geometric results about
SLE$(\kappa;\vec{\rho})$ processes. We will prove that many propositions in \cite{duality} and
Section \ref{GeometricProperties} of this paper about the limit of an SLE$(\kappa;\vec{\rho})$ trace that hold for $\kappa\in(0,4]$ are
also true for $\kappa>4$.

Julien Dub\'edat studied (Theorem 1, \cite{Julien-Duality}) the distribution of the boundary arc of $K(t)$ straddling $x$,
i.e., the boundary arc seen by $x$ at time $T_x^-$.  His result is about the ``inner'' boundary of $K(t)$, while Theorem
\ref{kappa(4,8)swallow} in this paper is about the ``outer'' boundary. The author feels that it is more appropriate
and convenient to apply Theorem  \ref{kappa(4,8)swallow}  to study the boundary of standard chordal SLE$(\kappa)$ hulls at 
general stopping times, and to derive other related results.

\section{Preliminary} \label{Preliminary}
In this section, we review some definitions and propositions in \cite{duality}, which will be used
in this paper.

If $H$ is a bounded and relatively  closed
subset of $\HH=\{z\in\C:\Imm z>0\}$, and $\HH\sem H$ is simply
connected, then we call $H$ a hull in $\HH$ w.r.t.\ $\infty$. For
such $H$, there is $\vphi_H$ that maps $\HH\sem H$ conformally onto
$\HH$, and satisfies $\vphi_H(z)=z+\frac{c}{z}+O(\frac 1{z^2})$ as
$z\to\infty$, where $c=\hcap(H)\ge 0$ is called the capacity of $H$
in $\HH$ w.r.t.\ $\infty$.

For a real interval $I$, we use $C(I)$ to denote the space of real
continuous functions on $I$. For $T>0$ and $\xi\in C([0,T))$, the
chordal Loewner equation driven by $\xi$ is
$$\pa_t\vphi(t,z)=\frac{2}{\vphi(t,z)-\xi(t)},\quad \vphi(0,z)=z.$$
For $0\le t<T$, let $K(t)$ be the set of $z\in\HH$ such that the
solution $\vphi(s,z)$ blows up before or at time $t$. We call $K(t)$
and $\vphi(t,\cdot)$, $0\le t<T$, chordal Loewner hulls and maps,
respectively, driven by $\xi$. It turns out that
$\vphi(t,\cdot)=\vphi_{K(t)}$ for each $t\in[0,T)$.

Let $B(t)$, $0\le t<\infty$, be a (standard linear) Brownian motion.
Let $\kappa \ge 0$. Then $K(t)$ and $\vphi(t,\cdot)$, $0\le
t<\infty$, driven by $\xi(t)=\sqrt\kappa B(t)$, $0\le t<\infty$, are
called standard chordal SLE$(\kappa)$ hulls and maps, respectively.
It is known (\cite{RS-basic}\cite{LSW-2}) that almost surely for any
$t\in[0,\infty)$, \BGE\gamma(t):=\lim_{\HH\ni
z\to\xi(t)}\vphi(t,\cdot)^{-1}(z)\label{trace}\EDE exists, and
$\gamma(t)$, $0\le t<\infty$, is a continuous curve in $\lin{\HH}$.
Moreover, if $\kappa\in(0,4]$ then $\gamma$ is a simple curve, which
intersects $\R$ only at the initial point, and for any $t\ge 0$,
$K(t)=\gamma((0,t])$; if $\kappa>4$ then $\gamma$ is not simple; if
$\kappa\ge 8$ then $\gamma$ is space-filling.
Such $\gamma$ is called a standard chordal SLE$(\kappa)$
trace.

If $(\xi(t))$ is a semi-martingale, and $d\langle \xi(t)\rangle
=\kappa dt$ for some $\kappa>0$, then from Girsanov theorem (c.f.\ \cite{RY}) and the
existence of standard chordal SLE$(\kappa)$ trace, almost surely for
any $t\in[0,T)$, $\gamma(t)$ defined by (\ref{trace}) exists, and has
the same property as a standard chordal SLE$(\kappa)$ trace
(depending on the value of $\kappa$) as described in the last
paragraph.

Let $\kappa\ge 0$, $\rho_1,\dots,\rho_N\in\R$, $x\in\R$, and
$p_1,\dots,p_N\in\ha\R\sem\{x\}$, where $\ha\R=\R\cup\{\infty\}$ is
a circle. Let $\xi(t)$ and $p_k(t)$, $1\le k\le N$, be the solutions
to the SDE:
\begin{equation}\left\{\begin{array}{lll} d\xi(t) & = &
\sqrt\kappa d B(t)+\sum_{k=1}^N\frac{\rho_k}{\xi(t)-p_k(t)}\,dt\\
dp_k(t) & = & \frac{2}{p_k(t)-\xi(t)}\,dt,\quad 1\le k\le
N,\end{array}\right.\label{kappa-rho}\end{equation} with initial
values $\xi(0)=x$ and $p_k(0)=p_k$, $1\le k\le N$. If
$\vphi(t,\cdot)$ are chordal Loewner maps driven by $\xi(t)$, then
$p_k(t)=\vphi(t,p_k)$. Suppose $[0,T)$ is the
maximal interval of the solution. Let $K(t)$ and $\gamma(t)$, $0\le t<T$, be chordal
Loewner hulls and trace driven by $\xi$. Let $\vec\rho=(\rho_1,\dots,\rho_N)$
and $\vec p=(p_1,\dots,p_N)$. Then  $K(t)$ and $\gamma(t)$, $0\le t<T$, are called
(full) chordal SLE$(\kappa;\rho_1,\dots,\rho_N)$ or SLE$(\kappa;\vec{\rho})$ process and trace, respectively, started
from $(x;p_1,\dots,p_N)$ or $(x;\vec{p})$.  If $T_0\in(0,T]$ is a stopping time, then
 $K(t)$ and $\gamma(t)$, $0\le t<T_0$, are called partial chordal
SLE$(\kappa;\vec{\rho})$ process and trace, respectively, started
from $(x;\vec{p})$.

If we allow that one of the force points takes value $x^+$ or $x^-$, or two of the
force points take values $x^+$ and $x^-$, respectively, then we obtain
the definition of degenerate chordal SLE$(\kappa;\vec{\rho})$ process.
Let $\kappa \ge 0$;
$\rho_1,\dots,\rho_N\in\R$, and $\rho_1\ge \kappa/2-2$; $p_1=x^+$,
$p_2,\dots,p_N\in \ha\R\sem\{x\}$. Let $\xi(t)$ and $p_k(t)$, $1\le
k\le N$, $0<t<T$, be the maximal solution to (\ref{kappa-rho}) with
initial values $\xi(0)=p_1(0)=x$, and $p_k(0)=p_k$, $1\le k\le N$.
Moreover, we require that $p_1(t)>\xi(t)$ for any $0<t<T$.
Then the chordal Loewner hulls $K(t)$ and trace $\gamma(t)$, $0\le t<T$, driven by $\xi$, are
called chordal  SLE$(\kappa;\rho_1,\dots,\rho_N)$ process and trace started from
$(x;x^+,p_2,\dots,p_N)$. If the condition $p_1(t)>\xi(t)$ is
replaced by $p_1(t)<\xi(t)$, then we get chordal
SLE$(\kappa;\rho_1,\dots,\rho_N)$ process and trace started from
$(x;x^-,p_2,\dots,p_N)$. Now suppose $N\ge 2$, $\rho_1,\rho_2\ge
\kappa/2-2$, $p_1=x^+$, and $p_2=x^-$. Let $\xi(t)$ and $p_k(t)$,
$1\le k\le N$, $0<t<T$, be the maximal solution to (\ref{kappa-rho})
with initial values $\xi(0)=p_1(0)=p_2(0)=x$, and $p_k(0)=p_k$,
$1\le k\le N$, such that $p_1(t)>\xi(t)>p_2(t)$ for all $0<t<T$.
Then we obtain chordal SLE$(\kappa;\rho_1,\dots,\rho_N)$ process and
trace started from $(x;x^+,x^-,p_3,\dots,p_N)$.
The force point $x^+$ or $x^-$ is called a degenerate force point.
Other force points are called generic force points. Let
$\vphi(t,\cdot)$ be the chordal Loewner maps driven by $\xi$. Since
for any generic force point $p_j$, we have $p_j(t)=\vphi(t,p_j)$, so
we write $\vphi(t,p_j)$ for $p_j(t)$ in the case
that $p_j$ is a degenerate force point.

For $h>0$, let
$\St_h=\{z\in\C:0<\Imm z<h\}$ and $\R_h=ih+\R$. If $H$ is a bounded
closed subset of $\St_\pi$,  $\St_\pi\sem H$ is simply connected,
and has $\R_\pi$ as a boundary arc, then we call $H$ a hull in
$\St_\pi$ w.r.t.\ $\R_\pi$. For such $H$, there is a unique $\psi_H$
that maps $\St_\pi\sem H$ conformally onto $\St_\pi$, such that for
some $c\ge 0$, $\psi_H(z)=z\pm c +o(1)$ as $z\to\pm\infty$ in
$\St_\pi$. We call such $c$ the capacity of $H$ in $\St_\pi$ w.r.t.\
$\R_\pi$, and let it be denoted it by $\scap(H)$.

For $\xi\in C([0,T))$, the strip Loewner equation driven by $\xi$ is
\BGE \pa_t\psi(t,z)=\coth\big(\frac{\psi(t,z)-\xi(t)}2\big),\quad
\psi(0,z)=z.\label{strip}\EDE For $0\le t<T$, let $L(t)$ be the set
of $z\in\St_\pi$ such that the solution $\psi(s,z)$ blows up before
or at time $t$. We call $L(t)$ and $\psi(t,\cdot)$, $0\le t<T$,
strip Loewner hulls and maps, respectively, driven by $\xi$. It
turns out that $\psi(t,\cdot)=\psi_{L(t)}$ and $\scap(L(t))=t$ for
each $t\in[0,T)$. In this paper, we use  $\coth_2(z)$, $\tanh_2(z)$,
$\cosh_2(z)$, and $\sinh_2(z)$ to denote the functions $\coth(z/2)$, $\tanh(z/2)$,
$\cosh(z/2)$, and $\sinh(z/2)$, respectively.

Let $\kappa\ge 0$, $\rho_1,\dots,\rho_N\in\R$, $x\in\R$, and
$p_1,\dots,p_N\in \R\cup\R_\pi\cup\{+\infty,-\infty\}\sem\{x\}$. Let
$B(t)$ be a Brownian motion. Let $\xi(t)$ and $p_k(t)$, $1\le k\le
N$, be the solutions to the SDE:
\begin{equation}\left\{\begin{array}{lll} d\xi(t) & = &
\sqrt\kappa d B(t)+\sum_{k=1}^N\frac{\rho_k}2\coth_2({\xi(t)-p_k(t)})dt\\
dp_k(t) & = & \coth_2({p_k(t)-\xi(t)})dt,\quad 1\le k\le
N,\end{array}\right.\label{strip-kappa-rho}\end{equation} with
initial values $\xi(0)=x$ and $p_k(0)=p_k$, $1\le k\le N$. Here if
some $p_k=\pm\infty$ then $p_k(t)=\pm\infty$ and
$\coth_2({\xi(t)-p_k(t)})=\mp 1$ for all $t\ge 0$. Suppose $[0,T)$ is the maximal
interval of the solution to (\ref{strip-kappa-rho}). Let $L(t)$, $0\le t<T$, be strip Loewner
hulls driven by $\xi$. Let $\beta(t)=\lim_{\St_\pi\ni z\to\xi(t)}\psi(t,z)$,
$0\le t<T$. Then we call $L(t)$ and $\beta(t)$, $0\le t<T$, (full)
strip SLE$(\kappa;\vec{\rho})$ process and trace, respectively, started from $(x;\vec{p})$,
where $\vec{\rho}=(\rho_1,\dots,\rho_N)$ and
$\vec{p}=(p_1,\dots,p_N)$. If $T_0\in(0,T]$ is a stopping time, then
 $L(t)$ and $\beta(t)$, $0\le t<T_0$, are called partial strip
SLE$(\kappa;\vec{\rho})$ process and trace, respectively, started
from $(x;\vec{p})$.

The following two propositions are Lemma 2.1 and Lemma 2.3 in \cite{duality}.
They will be used frequently in this paper. Let $S_1$ and $S_2$ be two sets of
boundary points or prime ends of a domain $D$. We say that $K$ does not
separate $S_1$ from $S_2$ in $D$ if there are neighborhoods $U_1$
and $U_2$ of $S_1$ and $S_2$, respectively, in $D$ such that $U_1$
and $U_2$ lie in the same pathwise connected component of $D\sem K$.

\begin{Proposition} Suppose $\kappa\ge 0$ and
$\vec{\rho}=(\rho_1,\dots,\rho_N)$ with  $\sum_{m=1}^N\rho_m =\kappa-6$.
For $j=1,2$, let $K_j(t)$, $0\le
t<T_j$, be a generic or degenerate chordal SLE$(\kappa;\vec{\rho})$ process started from
$(x_{j};\vec{p}_{j})$, where $\vec{p}_{j}=(p_{j,1},\dots,p_{j,N})$,
$j=1,2$. Suppose $W$ is a conformal
or conjugate conformal map from $\HH$ onto $\HH$ such that $W(x_1)=x_2$ and
$W(p_{1,m})=p_{2,m}$, $1\le m\le N$. Let
$p_{1,\infty}=W^{-1}(\infty)$ and $p_{2,\infty}=W(\infty)$. For
$j=1,2$, let $S_j\in(0,T_j]$ be the largest number such that for $0\le
t<S_j$,  $K_j(t)$ does not separate $p_{j,\infty}$ from $\infty$ in
$\HH$. Then $(W(K_1(t)),0\le t<S_1)$ has the same law as
$(K_2(t),0\le t<S_2)$ up to a time-change. A similar result holds
for the traces.  \label{coordinate}
\end{Proposition}

\begin{Proposition} Suppose $\kappa\ge 0$ and $\vec{\rho}=(\rho_1,\dots,\rho_N)$
with $\sum_{m=1}^N\rho_m =\kappa-6$.   Let $K(t)$, $0\le t<T$, be a
chordal SLE$(\kappa;\vec{\rho})$ process started from
$(x;\vec{p})$, where $\vec{p}=(p_{1},\dots,p_{N})$. Let $L(t)$, $0\le t<S$, be a strip
SLE$(\kappa;\vec{\rho})$ process started from $(y;\vec{q})$, where
$\vec{q}=(q_{1},\dots,q_{N})$. Suppose $W$ is a conformal
or conjugate conformal map from $\HH$ onto $\St_\pi$ such that $W(x)=y$ and
$W(p_k)=q_k$, $1\le k\le N$. Let $I=W^{-1}(\R_\pi)$ and
$q_\infty=W(\infty)$. Let $T'\in(0,T]$ be the largest number such that for
$0\le t<T'$, $K(t)$ does not separate $I$ from $\infty$ in $\HH$.
Let $S'\in(0,S]$ be the largest number such that for $0\le t<S'$, $L(t)$
does not separate $q_\infty$ from $\R_\pi$. Then $(W(K(t)),0\le
t<T')$ has the same law as $(L(t),0\le t<S')$ up to a time-change.
A similar result holds for the traces.
\label{coordinate2}
\end{Proposition}

Now we recall some geometric results of SLE$(\kappa;\vec{\rho})$ traces
derived in \cite{duality}.

Let $\kappa>0$, and $\rho_+,\rho_-\in\R$ be such that
$\rho_++\rho_-=\kappa-6$. Suppose $\beta(t)$, $0\le t<\infty$, is a
strip SLE$(\kappa;\rho_+,\rho_-)$ trace
started from $(0;+\infty,-\infty)$. In the following propositions, Proposition \ref{less2} is
a combination of Lemma 3.1, Lemma 3.2, and the argument before Lemma 3.2, in \cite{duality};
Proposition \ref{less2'} and Proposition \ref{greater2} are
Theorem 3.3, and Theorem 3.4, respectively, in \cite{duality}.

\begin{Proposition} If $|\rho_+-\rho_-|<2$, then a.s.\ $\beta([0,\infty))$ is bounded, and
$\lin{\beta([0,\infty))}$ intersects $\R_\pi$
at a single point $J+\pi i$. And the distribution of $J$ has a probability density function w.r.t.\
the Lebesgue measure, which
is proportional to $\exp(x/2)^{\frac{2}\kappa\,(\rho_--\rho_+)} \cosh_2(x)^{-\frac 4\kappa}$.
\label{less2}
\end{Proposition}

\begin{Proposition} If $\kappa\in(0,4]$ and $|\rho_+-\rho_-|<2$,  then a.s.\
$\lim_{t\to \infty} \beta(t)\in\R_\pi$. \label{less2'}
\end{Proposition}

\begin{Proposition} If $\kappa\in(0,4]$ and $\pm(\rho_+-\rho_-)\ge 2$, then a.s.\
$\lim_{t\to \infty} \beta(t)=\mp\infty$.\label{greater2}
\end{Proposition}

The following two propositions are Theorem 3.1 and Theorem 3.2 in \cite{duality}.

\begin{Proposition} Let $\kappa>0$, $N_+,N_-\in\N$,
$\vec{\rho}_\pm=(\rho_{\pm 1},\dots,\rho_{\pm N_\pm})\in\R^{N_\pm}$ with
$\sum_{j=1}^k\rho_{\pm j}\ge \kappa/2-2$ for $1\le k\le N_\pm$,
$\vec{p}_\pm=(p_{\pm 1},\dots,p_{\pm N_\pm})$ with $0<p_1<\dots<p_{N_+}$ and
$0>p_{-1}>\dots>p_{-N-}$. Let
$\gamma(t)$, $0\le t<T$, be a chordal
SLE$(\kappa;\vec{\rho}_+,\vec{\rho}_-)$ trace started from
$(0;\vec{p}_+,\vec{p}_-)$.  Then a.s.\ $T=\infty$ and  $\infty$ is a
subsequential limit of $\gamma(t)$ as
$t\to \infty$.
\label{more force} \end{Proposition}

\begin{Proposition}
Let $\kappa\in(0,4]$, $\rho_+,\rho_-\ge \kappa/2-2$. Suppose $\gamma(t)$, $0\le t<\infty$,
is a chordal SLE$(\kappa;\rho_+,\rho_-)$ trace started from $(0;p_+,p_-)$. If $p_+ =0^+$ and
$p_-=0^-$, or $p^+\in(0,\infty)$ and $p^-\in(-\infty,0)$,
then a.s.\
$\lim_{t\to\infty}\gamma(t)=\infty$. \label{3pt**}
\end{Proposition}

The following proposition is
Theorem 4.1 in \cite{duality} in the case that $\kappa_1<4<\kappa_2$.

\begin{Proposition} Let $0<\kappa_1<4<\kappa_2$ be such that $\kappa_1\kappa_2=16$.
Let $x_1\ne x_2\in\R$. Let $N\in\N$. Let $p_1,\dots,p_N\in\R\cup\{\infty\}\sem\{x_1,x_2\}$ be distinct points.
For $1\le m\le N$, let $C_m\in\R$ and $\rho_{j,m}=C_m(\kappa_j-4)$,
$j=1,2$. There is a coupling of $K_1(t)$,
$0\le t<T_1$, and $K_2(t)$, $0\le t<T_2$, such that (i) for $j=1,2$,
$K_j(t)$, $0\le t<T_j$, is a chordal
SLE$(\kappa_j;-\frac{\kappa_j}2,\vec{\rho}_j)$ process started from
$(x_j;x_{3-j},\vec{p})$; and (ii) for $j\ne k\in\{1,2\}$, if $\bar
t_k$ is an $(\F^k_t)$-stopping time with $\bar t_k<T_k$, then
conditioned on $\F^k_{\bar t_k}$, $\vphi_k(\bar t_k,K_j(t))$, $0\le
t\le T_j(\bar t_k)$, has the same distribution as a time-change of a
partial chordal SLE$(\kappa_j;-\frac{\kappa_j}2,\vec{\rho}_j)$
process started from $(\vphi_k(\bar t_k,x_j);\xi_k(\bar
t_k),\vphi_k(\bar t_k,\vec{p}))$, where
$\vphi_k(t,\vec{p})=(\vphi_k(t,p_1),\dots,\vphi_k(t,p_N))$;
$\vphi_k(t,\cdot)$, $0\le t<T_k$, are chordal Loewner maps for the hulls
$K_k(t)$, $0\le t<T_k$; $T_j(\bar t_k)\in(0,T_j]$ is the
largest number such that $\lin{K_j(t)}\cap\lin{K_k(\bar t_k)}=\emptyset$
for $0\le t<T_j(\bar t_k)$; and $(\F^j_t)$ is the filtration
generated by $(K_j(t))$, $j=1,2$. This still holds if some $p_m$ take(s) value $x_1^\pm$ or $x_2^\pm$.
\label{bartk}
\end{Proposition}

\section{Integration of SLE measures} \label{SLE-intergration}
Let $\kappa>0$, $\rho_+,\rho_-\in\R$, $\rho_++\rho_-=\kappa-6$, and $|\rho_+-\rho_-|<2$.
Suppose $\xi(t)$, $0\le t<\infty$, is the driving function of
a strip SLE$(\kappa;\rho_+,\rho_-)$ process started from $(0;+\infty,-\infty)$.
Let $\sigma=(\rho_--\rho_+)/2$. Then there is a Brownian motion $B(t)$ such that
$\xi(t)=B(t)+\sigma t$, $0\le t<\infty$.

Let $\mu$ denote the distribution of $\xi$.
We consider $\mu$ as a probability measure on $C([0,\infty))$.
Let $(\F_t)$ be the filtration on $C([0,\infty))$ generated by coordinate maps.
Then the total $\sigma$-algebra is $\F_\infty=\vee_{t\ge 0}\F_t$.
For each $x\in\R$, let $\nu_x$ denote the
the distribution of the driving function of a strip SLE$(\kappa;-4,\rho_-+2,\rho_++2)$ process
started from $(0;x+\pi i,+\infty,-\infty)$, which is also a probability measure
on $C([0,\infty))$. Then we have the following lemma.

\begin{Lemma} We have
$$\mu=\frac{1}{Z}\int_{\R} \nu_x \exp(x/2)^{\frac{4}\kappa\,\sigma} \cosh_2(x)^{-\frac 4\kappa} dx,$$
where $dx$ is Lebesgue measure, $Z=\int_{\R}\exp(x/2)^{\frac{4}\kappa\,\sigma} \cosh_2(x)^{-\frac 4\kappa} dx$,
which is finite because $|\sigma|<1$, and the integral means that for any $A\in\F_\infty$,
\BGE \mu(A)=\frac{1}{Z}\int_{\R} \nu_x(A) \exp(x/2)^{\frac{4}\kappa\,\sigma}
 \cosh_2(x)^{-\frac 4\kappa}dx\label{integral2}.\EDE
\label{integration}
\end{Lemma}
{\bf Proof.} Let $f(x)=\frac 1Z\exp(x/2)^{\frac{4}\kappa\,\sigma} \cosh_2(x)^{-\frac 4\kappa}$, $x\in\R$.
Then $\int_\R f(x) dx=1$, and
\BGE \frac{f'(x)}{f(x)}=\frac{2}\kappa(\sigma-\tanh_2(x)),\quad x\in\R;\label{f'/f}\EDE
\BGE\frac{\kappa}2\,f''(x)+f'(x)(-\sigma+\tanh_2(x))+\frac{f(x)}2\,\cosh_2(x)^{-2}=0,
\quad x\in\R.\label{f}\EDE
Note that the collection of $A$ that satisfies (\ref{integral2}) is a monotone class,
and $\cup_{t\ge 0}\F_t$ is an algebra. From Monotone Class Theorem, we suffice to show that
 (\ref{integral2}) holds for any $A\in\F_t$, $t\in[0,\infty)$. This will be proved by showing
that $\nu_x|_{\F_t}\ll \mu|_{\F_t}$ for all $x\in\R$ and $t\in[0,\infty)$, and if
$R_t(x)$ is the Radon-Nikodym derivative, then $\int_\R R_t(x)f(x)dx=1$.

Let $\psi(t,\cdot)$, $0\le t<\infty$, be the strip Loewner maps driven by $\xi$. For $x\in\R$ and
$t\ge 0$, let $X(t,x)=\Ree(\psi(t,x+\pi i)-\xi(t))$. Note that $\psi(t,x+\pi i)\in\R_\pi$ for
any $t\ge 0$. From (\ref{strip}), for any fixed $x\in\R$, $X(t,x)$ satisfies the SDE
\BGE \pa_t X(t,x)=-\sqrt\kappa \pa B(t)-\sigma \pa t +\tanh_2(X(t,x))\pa t.\label{X-SDE}\EDE
If $t$ is fixed, then $\pa_x X(t,x)=\pa_x\psi(t,x+\pi i)$. From (\ref{strip}), we have
$$\pa_t\pa_x X(t,x)=\pa_t\pa_x\psi(t,x+\pi i)
=\frac 12\sinh_2(\psi(t,x+\pi i)-\xi(t))^{-2}\pa_x\psi(t,x+\pi i)$$
\BGE =\frac 12\cosh_2(X(t,x))^{-2}\pa_x X(t,x).\label{X'}\EDE
For $x\in\R$ and $t\ge 0$, define
$M(t,x)=f(X(t,x))\pa_x X(t,x)$.
From (\ref{f'/f}$\sim$\ref{X'}) and Ito's formula (c.f.\ \cite{RY}), we find that for any fixed $x$,
$(M(t,x))$ is a local martingale, and satisfies the SDE:
$$\frac{\pa_t M(t,x)}{M(t,x)}=-\frac{f'(X(t,x))}{f(X(t,x))}\,\sqrt\kappa \pa B(t)
=-\frac 2{\sqrt\kappa}(\sigma-\tanh_2(X(t,x)))\, \pa B(t).$$
From the definition, $f$ is bounded on $\R$. From (\ref{X'}) and that $\pa_x X(0,x)=1$, we have
$|\pa_x X(t,x)|\le \exp(t/2)$. Thus for any fixed $t_0>0$, $M(t,x)$ is bounded on $[0,t_0]\times\R$.
So $(M(t,x):0\le t\le t_0)$ is a bounded martingale.
Then we have $\EE[M(t_0,x)]=M(0,x)=f(x)$ for any $x\in\R$. Now define the probability measure $\nu_{t_0,x}$ such that
$d\nu_{t_0,x}/d\mu=M(t_0,x)/f(x)$, and let
$$\til B(t)=B(t)+\int_0^t \frac 2{\sqrt\kappa}(\sigma-\tanh_2(X(s,x)))\,ds, \quad
0\le t\le t_0.$$ From Girsanov Theorem, under the probability measure $\nu_{t_0,x}$,
$\til B(t)$, $0\le t\le t_0$, is a partial
Brownian motion. Now $\xi(t)$, $0\le t\le t_0$, satisfies the SDE:
$$d\xi(t)=\sqrt\kappa d\til B(t)+\sigma dt-2(\sigma-\tanh_2(X(t,x)))dt$$
$$=\sqrt\kappa d\til B(t)-\sigma dt-\frac{-4}2\,\coth_2(\psi(t,x+\pi i)-\xi(t))dt.$$
Since $\xi(0)=0$, so under $\nu_{t_0,x}$, $(\xi(t),0\le t\le t_0)$ has the distribution of the driving
function of a strip SLE$(\kappa;-4,\rho_-+2,\rho_++2)$ process started from $(0;x+\pi i,+\infty,-\infty)$.
So we conclude that $\nu_{t_0,x}|_{\F_{t_0}}=\nu_{x}|_{\F_{t_0}}$. Thus $\nu_x|_{\F_{t_0}}\ll \mu|_{\F_{t_0}}$,
and the Radon-Nikodym derivative is $R_{t_0}(x)=M(t_0,x)/f(x)$. Thus
$$\int_{\R} R_{t_0}(x)f(x)dx=\int_{\R} M(t_0,x)dx=\int_{\R} f(X(t_0,x))\pa_x X(t_0,x)dx
=\int_{\R} f(y)dy=1.\quad \Box$$

\begin{Theorem} Let $\kappa>0$, and $\rho_+,\rho_-\in\R$ satisfy $\rho_++\rho_-=\kappa-6$ and $|\rho_+-\rho_-|<2$.
Let $\bar\mu$ denote the distribution of a strip SLE$(\kappa;\rho_+,\rho_-)$ trace $\beta(t)$, $0\le t<\infty$,
started
from $(0;+\infty,-\infty)$. Let $\lambda$ denote the distribution of the intersection point of $\lin{\beta([0,\infty))}$ with $\R_\pi$.
For each $p\in\R_\pi$, let $\bar\nu_p$ denote the distribution of
a strip SLE$(\kappa;-4,\rho_-+2,\rho_++2)$ trace started from $(0;p,+\infty,-\infty)$. Then
$\bar\mu=\int_{\R_\pi} \bar\nu_p\, d\lambda(p)$. \label{Theorem-integration}
\end{Theorem}
{\bf Proof.} This follows from Proposition \ref{less2} and the above lemma. $\Box$

\vskip 3mm

\no{\bf Remark.} A special case of the above theorem is that
$\kappa=2$ and $\rho_+=\rho_-=-2$, so $\rho_++2=\rho_-+2=0$.
From \cite{LERW}, a strip SLE$(2;-2,-2)$
trace started from $(0;+\infty,-\infty)$ is a continuous LERW in $\St_\pi$ from $0$ to $\R_\pi$;
a strip SLE$(2;-4,0,0)$ trace started from $(0;p,+\infty,-\infty)$ is a continuous LERW in
$\St_\pi$ from $0$ to $p$; and the above theorem in this special case follows from the convergence
of discrete LERW to continuous LERW.

\begin{Corollary}
Let $\kappa>0$, $\rho\in(\kappa/2-4,\kappa/2-2)$, and $x\ne 0$. Let $\bar\mu$ denote the distribution
of a chordal SLE$(\kappa;\rho)$ trace $\gamma(t)$, $0\le t<T$, started from $(0;x)$.
Let $\lambda$ denote the distribution of the subsequential limit of $\gamma(t)$ on $\R$ as $t\to T$, which is a.s.\ unique.
For each $y\in\R$, let $\bar\nu_y$ denote the distribution of a chordal SLE$(\kappa;-4,\kappa-4-\rho)$ trace started from
$(0;y,x)$. Then $\bar\mu=\int_\R \bar\nu_y\,d\lambda(y)$.
\end{Corollary}
{\bf Proof.} This follows from the above theorem and Proposition \ref{coordinate2}. $\Box$

\begin{Corollary}
Let $\kappa\in(4,8)$ and $x\ne 0$. Let $\gamma(t)$, $0\le t<\infty$, be a standard chordal SLE$(\kappa)$
trace. Let $T_x$ be the first $t$ that $\gamma([0,t])$ disconnects $x$ from $\infty$ in $\HH$. Let
 $\bar\mu$ denote the distribution of $(\gamma(t),0\le t<T_x)$. Let $\lambda$ denote the distribution
of $\gamma(T_x)$. For each $y\in\R$, let
$\bar\nu_y$ denote the distribution of a chordal SLE$(\kappa;-4,\kappa-4)$ trace started from
$(0;y,x)$. Then
$\bar\mu=\int_\R \bar\nu_y\,d\lambda(y)$. \label{kappa(4,8)}
\end{Corollary}
{\bf Proof.} This is a special case of the above corollary because
$\gamma(t)$, $0\le t<T_x$, is a chordal SLE$(\kappa;0)$ trace started from $(0;x)$,
and $0\in(\kappa/2-4,\kappa/2-2)$.
 $\Box$

\section{Geometric Properties} \label{GeometricProperties}
In this section, we will improve some results derived in Section 3 of \cite{duality}. We first
derive a simple lemma.

\begin{Lemma} Suppose $\psi(t,\cdot)$, $0\le t<T$, are strip Loewner maps driven by $\xi$.
Suppose $\xi(0)<x_1<x_2$ or $\xi(0)>x_1>x_2$, and $\psi(t,x_1)$ and $\psi(t,x_2)$ are defined
for $0\le t<T$. Then for any $0\le t<T$,
$$\Big|\int_0^t\coth_2(\psi(s,x_1)-\xi(s))ds-\int_0^t\coth_2(\psi(s,x_2)-\xi(s))ds\Big|< |x_1-x_2|.$$
\label{close}
\end{Lemma}
{\bf Proof.} By symmetry, we only need to consider the case that $\xi(0)<x_1<x_2$. For any $0\le t<T$, we have
$\xi(t)<\psi(t,x_1)<\psi(t,x_2)$, which implies that $\coth_2(\psi(t,x_1)-\xi(t))>\coth_2(\psi(t,x_2)-\xi(t))>0$.
Also note that $\pa_t\psi(t,x_j)=\coth_2(\psi(t,x_j)-\xi(t))$, $j=1,2$, so for $0\le t<T$,
$$0\le \int_0^t\coth_2(\psi(s,x_1)-\xi(s))ds-\int_0^t\coth_2(\psi(s,x_2)-\xi(s))ds$$
$$=(\psi(t,x_1)-\psi(0,x_1))-(\psi(t,x_2)-\psi(0,x_2))$$$$=\psi(t,x_1)-\psi(t,x_2)+x_2-x_1<x_2-x_1
=|x_1-x_2|.\quad\Box$$

\vskip 3mm
From now on, in this section, we let $\kappa>0$, $N_+,N_-\in\N\cup\{0\}$,
$\vec{\rho}_\pm=(\rho_{\pm 1},\dots,\rho_{\pm N_\pm})\in\R^{N_\pm}$, and
$\chi_\pm =\sum_{m=1}^{N_\pm}\rho_{\pm m}$. Let $\tau_+,\tau_-\in\R$ be such that
$\chi_++\tau_++\chi_-+\tau_-=\kappa-6$. Let $\vec{p}_\pm=(p_{\pm 1},\dots,p_{\pm N_\pm})$
be such that $p_{-N_-}<\dots<p_{-1}<0<p_1<\dots<p_{N_+}$. Suppose $\beta(t)$, $0\le t<T$,
is a strip SLE$(\kappa;\vec{p}_+,\vec{p}_-,\tau_+,\tau_-)$ trace started from $(0;
\vec{p}_+,\vec{p}_-,+\infty,-\infty)$. Let $\xi(t)$ and $\psi(t,\cdot)$, $0\le t<T$, be the driving
function and strip Loewner maps for $\beta$. Then there is a Brownian motion $B(t)$ such
that for $0\le t<T$, $\xi(t)$ satisfies the SDE
$$d\xi(t)=\sqrt\kappa dB(t)-\sum_{m=1}^{N_+}\frac{\rho_m}2\coth_2(\psi(t,p_m)-\xi(t))\,dt$$
\BGE -\sum_{m=1}^{N_-}\frac{\rho_{-m}}2\coth_2(\psi(t,p_{-m})-\xi(t))\,dt-\frac{\tau_+-\tau_-}2\,dt.
\label{SDE}\EDE
For $0\le t<T$, we have
\BGE \psi(t,p_{-N_-})<\dots<\psi(t,p_{-1})<\xi(t)<\psi(t,p_1)<\dots<\psi(t,p_{N_+}).\label{order}\EDE
Since $\pa_t\psi(t,x)=\coth_2(\psi(t,x)-\xi(t))$, so $\pa_t\psi(t,p_m)>1$ for $1\le m\le N_+$,
and $\pa_t\psi(t,p_{-m})<-1$ for $1\le m\le N_-$. Thus for $0\le t<T$, $\psi(t,p_m)$ increases
in $t$, and $\psi(t,p_m)>t$ for $1\le m\le N_+$; $\psi(t,p_{-m})$ decreases in $t$, and
$\psi(t,p_{-m})<-t$ for $1\le m\le N_-$. We say that some force point
$p_s$ is swallowed by $\beta$  if $T<\infty$ and $\psi(t,p_s)-\xi(t)\to 0$ as $t\to T$.
In fact, if $T<\infty$ then some force point on $\R$ must be swallowed by $\beta$,
and from (\ref{order}) we see that either $p_1$ or $p_{-1}$ is swallowed.

\begin{Lemma} (i) If $\sum_{j=1}^k\rho_{j}\ge \kappa/2-2$ for $1\le k\le N_+$, then a.s.\ $p_1$ is not
swallowed by $\beta$. (ii) If $\sum_{j=1}^k\rho_{-j}\ge \kappa/2-2$ for $1\le k\le N_-$,
then a.s.\ $p_{-1}$ is not swallowed by $\beta$. \label{T=infty}
\end{Lemma}
{\bf Proof.} From symmetry we only need to prove (i). Suppose $\sum_{j=1}^k\rho_{j}\ge \kappa/2-2$
for $1\le k\le N_+$. Let $\cal E$ denote the event that $p_1$ is swallowed by $\beta$.
Let $\PP$ be the probability measure we are working on. We want to show that $\PP({\cal E})=0$.
Assume that $\PP({\cal E})>0$.
 Assume that $\cal E$ occurs. Then $\lim_{t\to T}\xi(t)=\lim_{t\to T}\psi(t,p_1)\ge T$. For $1\le m\le N_-$,
since $\psi(t,p_{-m})<-t$, $0\le t<T$, so $\psi(t,p_{-m})-\xi(t)$ on $[0,T)$ is uniformly bounded
above by a negative number. Thus $\coth_2(\psi(t,p_{-m})-\xi(t))$ on $[0,T)$ is uniformly bounded for
$1\le m\le N_-$. For $0\le t<T$, let $\til B(t)=B(t)+\int_0^t a(s)ds$, where
$$a(t)= -\frac{\kappa/2-2}{2\sqrt\kappa}+\frac{\kappa/2-4-\chi_+}{2\sqrt\kappa}\,\coth_2(\psi(t,\pi i)-\xi(t))$$ $$-\sum_{m=1}^{N_-}\frac{\rho_{-m}}{2\sqrt\kappa}\coth_2(\psi(t,p_{-m})-\xi(t))-\frac{\tau_+-\tau_-}
{2\sqrt\kappa}.$$
For $0\le t<T$, since $\psi(t,\pi i)-\xi(t)\in\R_\pi$, so $|\coth_2(\psi(t,\pi i)-\xi(t))|\le 1$.
From the previous discussion, we see that if $\cal E$ occurs, then $T<\infty$ and $a(t)$ is uniformly bounded
on $[0,T)$, and so $\int_0^T a(t)^2 dt<\infty$. For $0\le t<T$, define
\BGE M(t)=\exp\Big(-\int_0^t a(s)dB(s)-\int_0^t a(s)^2ds\Big).\label{M(t)}\EDE
Then $(M(t),0\le t<T)$ is a local martingale and satisfies $dM(t)/M(t)=-a(t)dB(t)$. In the
event $\cal E$, since $\int_0^T a(t)^2 dt<\infty$, so a.s.\ $\lim_{t\to T} M(t)\in (0,\infty)$.
For $N\in\N$,  let $T_N\in[0,T]$ be the largest number such that $M(t)\in(1/{(2N)},2N)$ on
$[0,T_N)$. Let ${\cal E}_N={\cal E}\cap\{T_N=T\}$. Then ${\cal E}=\cup_{N=1}^\infty {\cal E}_N$ a.s.,
and $\EE[M(T_N)]=M(0)=1$, where $M(T):=\lim_{t\to T}M(t)$.  Since $\PP({\cal E})>0$, so there is
$N\in\N$ such that $\PP({\cal E}_N)>0$. Define another probability measure $\QQ$ such that
$d\QQ/d\PP=M(T_N)$. Then $\PP\ll \QQ$, and so
$\QQ({\cal E}_N)>0$. By Girsanov Theorem,
under the probability measure $\QQ$, $\til B(t)$, $0\le t<T_N$, is a partial Brownian motion.
From (\ref{SDE}), $\xi(t)$, $0\le t<T$, satisfies the SDE:
$$d\xi(t)=\sqrt\kappa d\til B(t)-\sum_{m=1}^{N_+}\frac{\rho_m}2\coth_2(\psi(t,p_m)-\xi(t))\,dt$$
$$+\frac{\kappa/2-2}{2}\,dt-\frac{\kappa/2-4-\chi_+}{2}\,\coth_2(\psi(t,\pi i)-\xi(t))\,dt,$$
so under $\QQ$, $\beta(t)$, $0\le t<T_N$, is a partial
strip SLE$(\kappa;\vec{\rho_+},\frac{\kappa}2-2,\frac{\kappa}2-4-\chi_+)$ trace started from
$(0;\vec{p}_+,-\infty,\pi i)$. In the event ${\cal E}_N$, since $\psi(t,p_1)-\xi(t)\to 0$
as $t\to T_N=T$, so $\beta(t)$, $0\le t<T_N$, is a full trace under $\QQ$.
Note that $$\sum_{m=1}^{N_+}\rho_m+\Big(\frac{\kappa}2-2\Big)+
\Big(\frac{\kappa}2-4-\chi_+\Big)=\kappa-6.$$
From Proposition \ref{coordinate2}, Proposition
\ref{more force}, and that $\sum_{j=1}^k\rho_{j}\ge \kappa/2-2$
for $1\le k\le N_+$, we see that on ${\cal E}_N$,
 $\QQ$-a.s.\ $\pi i$ is a subsequential limit of $\beta(t)$ as $t\to T_N$,
which implies that the height of $\beta((0,t])$ tends to $\pi$ as $t\to T_N$, and so $T_N=\infty$.
This contradicts that $T_N= T<\infty$ on ${\cal E}_N$ and $\QQ({\cal E}_N)>0$. Thus $\PP({\cal E})=0$.
$\Box$

\begin{Lemma} (i) If $\chi_+\ge -2$ and $\chi_++\tau_+>\kappa/2-4$, then $T=\infty$ a.s.\ implies that
$\liminf_{t\to\infty}(\psi(t,p_{m})-\xi(t))/t>0$ for $1\le m\le N_+$. (ii) If $\chi_-\ge -2$ and
$\chi_-+\tau_->\kappa/2-4$, then a.s.\ $T=\infty$ implies that
$\limsup_{t\to \infty} (\psi(t,p_{-m})-\xi(t))/t<0$ for $1\le m\le N_-$. \label{X/t}
\end{Lemma}
{\bf Proof.} We will only prove (i) since (ii)
follows from symmetry. Suppose $\chi_+\ge -2$ and $\chi_+ +\tau_+ >\kappa/2-4$.
Then $\Delta:=1+\frac{\chi_+}2+\frac{\tau_+}2-\frac{\chi_-}2-\frac{\tau_-}2>0$. Let
$X(t)=\psi(t,p_1)-\xi(t)$, $t\ge 0$. From (\ref{order}) we suffice to
show that $T=\infty$ a.s.\ implies that $\liminf_{t\to \infty}X(t)/t>0$.
Now assume that $T=\infty$. From (\ref{strip}) and (\ref{SDE}), for any $0\le t_1\le t_2$,
$$ X(t_2)-X(t_1)=-\sqrt\kappa B(t_2)+\sqrt\kappa B(t_1)+\frac{\tau_+-\tau_-}2\,(t_2-t_1)+\int_{t_1}^{t_2}\coth_2(X(t))dt$$
\BGE +\sum_{m=1}^{N_+}\frac{\rho_m}2\int_{t_1}^{t_2}\coth_2(\psi(t,p_m)-\xi(t))dt
+\sum_{m=1}^{N_-}\frac{\rho_{-m}}2\int_{t_1}^{t_2}\coth_2(\psi(t,p_{-m})-\xi(t))dt.\label{X}\EDE
Let $M_+=\sum_{m=1}^{N_+}|\rho_m||p_m-p_1|$.
From Lemma \ref{close}, for any $0\le t_1\le t_2$,
\BGE \sum_{m=1}^{N_+}\frac{\rho_m}2\int_{t_1}^{t_2}\coth_2(\psi(t,p_m)-\xi(t))dt
\ge \frac{\chi_+}2 \int_{t_1}^{t_2}\coth_2(X(t))dt-M_+.\label{M+}\EDE
Let $\eps_1=\min\{\Delta,1\}/6>0$.  There is a random number $A_0=A_0(\omega)>0$ such that
a.s.
\BGE |\sqrt\kappa B(t)|\le A_0+\eps_1 t,\quad \mbox{for any }  t\ge 0.\label{B}\EDE
Let $\chi^*_-=\sum_{m=1}^{N_-}|\rho_{-m}|$, and $\eps_2=\frac{\Delta}{\chi^*_-+1}>0$.
Choose $R>0$ such that if $x<-R$ then $|\coth_2(x)-(-1)|<\eps_2$.
Suppose $X(t)\le t$ on $[t_1,t_2]$, where $t_2\ge t_1\ge R$. Then for $1\le m\le N_-$ and $t\in[t_1,t_2]$,
from $\psi(t,p_{-m})<-t$ and $\psi(t,p_1)>t$, we have
$$\psi(t,p_{-m})-\xi(t)=\psi(t,p_{-m})-\psi(t,p_1)+X(t)<-t-t+t=-t\le -R,$$
and so $|\coth_2(\psi(t,p_{-m})-\xi(t))-(-1)|<\eps_2$.  Then
\BGE\sum_{m=1}^{N_-}\frac{\rho_{-m}}2\int_{t_1}^{t_2}\coth_2(\psi(t,p_{-m})-\xi(t))dt\ge
\Big(-\frac{\chi_-}2-\frac{\chi_-^*}2\,\eps_2\Big)\,(t_2-t_1).\label{M-}\EDE

Suppose $X(t_0)\ge t_0$ for some $t_0\ge \max\{R,2M_++4A_0+2\}$. We claim that a.s.\ for any $t\ge t_0$, we
have $X(t)\ge \eps_1 t$. If this is not true, then there are $t_2>t_1\ge t_0$ such that $X(t_1)=t_1$,
$X(t_2)=\eps_1 t_2$ and $X(t)\le t$ for $t\in[t_1,t_2]$.
From (\ref{X}$\sim$\ref{M-}), we have a.s.\
$$X(t_2)-X(t_1)\ge -2A_0-\eps_1t_1-\eps_1 t_2+\frac{\tau_+-\tau_-}2\,(t_2-t_1)$$
$$+\Big(1+\frac{\chi_+}2\Big)\int_{t_1}^{t_2}\coth_2(X(t))dt
 -M_+-\frac{\chi_-+\chi_-^*\eps_2}2\,(t_2-t_1).$$
\BGE \ge -M_+ -2A_0-2\eps_1 t_2+(\Delta-\frac{\chi_-^*\eps_2}2)(t_2-t_1),\label{X-X}\EDE
where in the last inequality we use the facts that $\coth_2(X(t))>1$ and
$1+\frac{\chi_+}2\ge  0$. Since $X(t_1)=t_1$ and $X(t_2)=\eps_1 t_2$, so
we have
$$M_+ +2A_0\ge (\Delta-{\chi_-^*\eps_2}/2-3\eps_1)(t_2-t_1)+(1-3\eps_1)t_1.$$
Since $\Delta-{\chi_-^*\eps_2}/2-3\eps_1\ge \Delta-\Delta/2-\Delta/2\ge 0$ and $1-3\eps_1\ge 1/2$, so
$$M_+ +2A_0\ge t_1/2\ge t_0/2\ge (2M_+ +4A_0+2)/2=M_++2A_0+1,$$
which is a contradiction. Thus if $X(t_0)\ge t_0$ for some $t_0\ge \max\{R,2M_++4A_0+2\}$, then a.s.\
$X(t)\ge \eps_1 t$ for any $t\ge t_0$, and so $\liminf_{t\to\infty}X(t)/t\ge \eps_1>0$.
The other possibility is that $X(t_0)<t_0$ for all $t_0\ge \max\{R,2M_++4A_0+2\}$.
Let $t_1=\max\{R,2M_++4A_0+2\}$ and $t_2\ge t_1$. Then (\ref{X}$\sim$\ref{M-}) still hold, so we have
(\ref{X-X}) again. Let both sides of (\ref{X-X}) be divided by $t_2$ and let $t_2=t\to\infty$. Then we have
a.s.\
$$\liminf_{t\to\infty} X(t)/t\ge \Delta-{\chi_-^*\eps_2}/2-2\eps_1\ge \Delta/6>0.\quad \Box$$

\vskip 3mm

The following theorem improves Theorem 3.6 in \cite{duality}.

\begin{Theorem} If $\kappa\in(0,4]$, $\sum_{j=1}^k\rho_{\pm j}\ge \kappa/2-2$, $1\le k\le N_\pm$,
and $|\chi_++\tau_+-\chi_--\tau_-|<2$,
then a.s.\ $T=\infty$ and $\lim_{t\to \infty}\beta(t)\in\R_\pi$.\label{lim-1p3}
\end{Theorem}
{\bf Proof.} From Lemma \ref{T=infty}, a.s.\ neither $p_1$ nor $p_{-1}$ is
swallowed by $\beta$, so $T=\infty$. Since $|\chi_++\tau_+-\chi_--\tau_-|<2$ and
$\chi_++\tau_++\chi_-+\tau_-=\kappa-6$, so $\chi_\pm+\tau_\pm>\kappa/2-4$.
If $N_+\ge 1$, then $\chi_+=\sum_{m=1}^{N_+}\rho_{\pm}\ge \kappa/2-2\ge -2$, so
from Lemma \ref{X/t}, a.s.\ $\liminf_{t\to\infty}(\psi(t,p_{m})-\xi(t))/t>0$ for $1\le m\le N_+$.
If $N_+=0$, this is also true since there is nothing to check. Similarly,
$\limsup_{t\to \infty} (\psi(t,p_{-m})-\xi(t))/t<0$ for $1\le m\le N_-$.
For $0\le t<\infty$, let $\til B(t)=B(t)+\int_0^t a(s)ds$, where
$$a(t)=\sum_{m=1}^{N_+}\frac{\rho_m}{2\sqrt\kappa}(1-\coth_2(\psi(t,p_m)-\xi(t)))
-\sum_{m=1}^{N_-}\frac{\rho_{-m}}{2\sqrt\kappa}(1+\coth_2(\psi(t,p_{-m})-\xi(t))).$$
Then $\int_0^\infty a(t)^2 dt<\infty$, and $\xi(t)$, $0\le t<\infty$, satisfies the SDE:
\BGE d\xi(t)=\sqrt\kappa d\til B(t)-\frac{\tau_++\chi_+-\tau_--\chi_-}2\,dt.\label{SDE2}\EDE
For $0\le t<\infty$, define $M(t)$ by (\ref{M(t)}). Then $(M(t))$ is a local martingale,
satisfies the SDE: $dM(t)/M(t)=-a(t)dB(t)$, and a.s.\ $M(\infty):=\lim_{t\to\infty}M(t)\in(0,\infty)$. For
$N\in\N$, let $T_N\in[0,\infty]$ be the largest number such that $M(t)\in(1/{(2N)},2N)$ on
$[0,T_N)$. Then $\EE[M(T_N)]=M(0)=1$. Let ${\cal E}_N=\{T_N=\infty\}$. Let $\PP$ be the probability
measure we are working on. Fix $\eps>0$. There is $N\in\N$ such that $\PP[{\cal E}_N]>1-\eps$. Define
another probability measure $\QQ$ such that $d\QQ/d\PP=M(T_N)$. By Girsanov Theorem,
under $\QQ$, $\til B(t)$, $0\le t<T_N$, is a partial Brownian motion, which together with
(\ref{SDE2}) implies that $\beta(t)$, $0\le t<T_N$, is a partial strip SLE$(\kappa;\rho_+,\rho_-)$
trace started from $(0;+\infty,-\infty)$, where $\rho_\pm=\chi_\pm+\tau_\pm$.
Since $\rho_++\rho_-=\kappa-6$ and $|\rho_+-\rho_-|<2$, so from Proposition \ref{less2'}, $\QQ$-a.s.\
$\lim_{t\to T_N}\beta(t)\in\R_\pi$ on $\{T_N=\infty\}={\cal E}_N$. Since $\PP\ll\QQ$, so
($\PP$-)a.s.\ $\lim_{t\to T_N}\beta(t)\in\R_\pi$ on ${\cal E}_N$. Since $\PP[{\cal E}_N]>1-\eps$,
so the probability that $\lim_{t\to \infty}\beta(t)\in\R_\pi$ is greater than $1-\eps$. Since
$\eps>0$ is arbitrary, so ($\PP$-)a.s.\ $\lim_{t\to \infty}\beta(t)\in\R_\pi$. $\Box$

\vskip 3mm

The following Theorem improves Theorem 3.1 in \cite{duality} when $\kappa\in(0,4]$.

\begin{Theorem}
Suppose $\kappa\in(0,4]$; $N_+,N_-\in\N\cup\{0\}$; $\vec{\rho}_\pm=(\rho_{\pm 1},\dots,\rho_{\pm N_\pm})
\in \R^{N_\pm}$; $\sum_{j=1}^k\rho_{\pm j}\ge \kappa/2-2$, $1\le k\le N_\pm$; $\vec{p}_\pm=(p_{\pm 1},
\dots,p_{\pm N_\pm})\in \R^{N_\pm}$; $p_{-N_-}<\dots<p_{-1}<0<p_1<\dots<p_{N_+}$. Let $\gamma(t)$, $0\le t<T$, be a
chordal SLE$(\kappa;\vec{\rho}_+,\vec{\rho}_-)$ trace started from $(0;\vec{p}_+,\vec{p}_-)$. Then
a.s.\ $\lim_{t\to T}\gamma(t)=\infty$.\label{tendstoinfty}
\end{Theorem}
{\bf Proof.} If $N_+=N_-=0$ then $\gamma$ is a standard chordal SLE$(\kappa)$ trace, so the conclusion
follows from Theorem 7.1 in \cite{RS-basic}. If $N_+=0$ and $N_-=1$, or $N_+=1$ and $N_-=0$, the conclusion follows from Proposition \ref{coordinate2}
and Proposition \ref{greater2}. If $N_+=N_-=1$, this follows from
Proposition \ref{3pt**}. For other cases, we will prove the theorem by reducing the number of force points.

Now consider the case that $N_-=0$ and $N_+\ge 2$. Choose $W$ that maps
$\HH$ conformally onto $\St_\pi$ such that $W(0)=0$, $W(\infty)=-\infty$, and $W(p_{N_+})=+\infty$.
Let $N_+'=N_+-1$; $\vec{q}=(q_1,\dots,q_{N_+'})$, where $q_m=W(p_m)$, $1\le m\le N_+'$. Then $0<q_1<\dots<q_{N_+'}$. Let $\vec{\rho}=(\rho_1,\dots,\rho_{N_+'})\in\R^{N_+'}$. Then
$\sum_{j=1}^k \rho_j\ge \kappa/2-2$ for $1\le k\le N_+'$. Let $\chi_+=\sum_{m=1}^{N_+'} \rho_m$. Then
$\chi_+\ge \kappa/2-2\ge -2$. Let $\tau_+=\rho_{N_+}$ and $\tau_-=\kappa-6-\chi_+-\tau_+$. Then
$\chi_++\tau_++\tau_-=\kappa-6$ and $\chi_++\tau_+=\sum_{m=1}^{N_+}\rho_m\ge \kappa/2-2>\kappa/2-4$.
From Proposition \ref{coordinate2}, a time-change of
$W\circ\gamma(t)$, $0\le t<T$, say $\beta(t)$, $0\le t<S$,
is a strip SLE$(\kappa;\tau_-,\tau_+,\vec{\rho})$ trace started from $(0;-\infty,+\infty,\vec{q})$.
Let $\xi(t)$ and $\psi(t,\cdot)$, $0\le t<S$, be the driving function and strip Loewner maps for $\beta$.
Then there is a Brownian motion $B(t)$
such that for $0\le t<S$, $\xi(t)$ satisfies the SDE
$$d\xi(t)=\sqrt\kappa dB(t)-\frac{\tau_+-\tau_-}2\,dt-\sum_{m=1}^{N_+'}\frac{\rho_m}2\coth_2
(\psi(t,q_m)-\xi(t))\,dt.$$
From Lemma \ref{T=infty} and Lemma \ref{X/t},  a.s.\ $S=\infty$ and
$\liminf_{t\to\infty} (\psi(t,q_m)-\xi(t))/t>0$ for $1\le m\le N_+'$. Let $\til B(t)=B(t)+\int_0^t a(s)ds$,
where
$$a(t)=\sum_{m=1}^{N_+'}\frac{\rho_m}{2\sqrt\kappa}(1-\coth_2(\psi(t,p_m)-\xi(t))).$$
Then $\int_0^\infty a(t)^2dt<\infty$. Now $\xi(t)$ satisfies the SDE
$$d\xi(t)=\sqrt\kappa d\til B(t)-\frac{\chi_++\tau_+-\tau_-}2\,dt.$$
Note that $(\chi_++\tau_+)-\tau_-\ge 2$.
We observe that if $\til B(t)$ is a Brownian motion, then $\beta$ is a strip
SLE$(\kappa;\chi_++\tau_+,\tau_-)$ trace started from $(0;+\infty,-\infty)$, and so from Proposition
\ref{greater2}, we have $\lim_{t\to\infty}\beta(t)=-\infty$.
Using the argument at the end of the proof of Theorem \ref{lim-1p3}, we conclude that a.s.\ $\lim_{t\to\infty}\beta(t)=-\infty$, and so $\lim_{t\to T}\gamma(t)
=W^{-1}(-\infty)=\infty$.

For the case $N_-=1$ and $N_+\ge 2$, we define $W$ and $\beta$ as in the above case, and
conclude that $\lim_{t\to\infty}\beta(t)=-\infty$ using the same argument as above except
that now we use Proposition \ref{coordinate2} and the conclusion of this theorem in the case $N_+=N_-=1$ to prove that a.s.\
$\lim_{t\to\infty}\beta(t)=-\infty$. So again we conclude that a.s.\
$\lim_{t\to T}\gamma(t) =\infty$. The cases that $N_+\in\{0,1\}$ and $N_-\ge 2$ are symmetric
to the above two cases. For the case that $N_+,N_-\ge 2$. we define $W$ and $\beta$ as in the
case that $N_-=0$ and $N_+\ge 2$, and conclude that a.s.\ $\lim_{t\to\infty}\beta(t)=-\infty$ using the
same argument as in that case except that now we use Proposition \ref{coordinate2} and the conclusion of this theorem in the case
$N_-\ge 2$ and $N_+=1$. So we also have a.s.\ $\lim_{t\to T}\gamma(t)=\infty$. $\Box$

\section{Duality} \label{Duality}
Let $\gamma$ be a simple curve in a simply connected domain $\Omega$.
We call $\gamma$ a crosscut in $\Omega$ if its two ends
approach to two different boundary points or prime ends of $\Omega$.
We call $\gamma$ a degenerate crosscut in $\Omega$ if its two ends approach to the same boundary point
or prime end of $\Omega$. We call $\gamma$ a semi-crosscut in $\Omega$ if its one end
approaches to some boundary point or prime end of $\Omega$, and the other end stays inside $\Omega$.
In the above definitions, if $\Omega=\HH$, and no end of $\gamma$ is $\infty$, then $\gamma$ is
called a crosscut, or degenerate crosscut, or semi-crosscut, respectively, in $\HH$ on $\R$.
For example, $e^{i\theta}$, $0<\theta<\pi$, is a crosscut in $\HH$ on $\R$; $e^{i\theta}$, $0<\theta\le \pi/2$,
is a semi-crosscut in $\HH$ on $\R$; $i+e^{i\theta}$, $-\pi/2<\theta<3\pi/2$, is a degenerate crosscut
in $\HH$ on $\R$.
If $\gamma$ is a crosscut in $\HH$ on $\R$, $\HH\sem\gamma$ has two connected components.
We use $D_\HH(\gamma)$ to denote the bounded component.

In Proposition \ref{bartk}, let  $N=4$; choose
$p_1<x_1<p_3<p_4<x_2<p_2$; choose $C_2,C_4\ge 1/2$, let $C_1=1-C_2$, $C_3=1/2-C_4$,
and $\rho_{j,m}=C_m(\kappa_j-4)$, $1\le m\le 4$, $j=1,2$.
Let $K_j(t)$, $0\le t<T_j$, $j=1,2$, be given by Proposition \ref{bartk}.
Let $\vphi_j(t,\cdot)$ and $\gamma_j(t)$, $0\le
t<T_j$, $j=1,2$, be the corresponding chordal Loewner maps and traces.

Since $\kappa_1\in(0,4)$, so $\gamma_1(t)$, $0\le t<T_j$, is a simple curve, and $\gamma_1(t)\in
\HH$ for $0<t<T_j$. From Theorem \ref{lim-1p3} and Proposition
\ref{coordinate2}, a.s.\
$\gamma_1(T_1):=\lim_{t\to T_1} \gamma_1(t)\in (x_2,p_2)$. Thus $\gamma_1$ is a
crosscut in $\HH$ on $\R$. Note that $\gamma_1$ disconnects $x_2$ from $\infty$ in $\HH$.
If $\bar t_2\in[0,T_2)$ is an $(\F^2_t)$-stopping time, then conditioned on $\F^2_{\bar t_2}$,
after a time-change, $\vphi_2(\bar t_2,\gamma_1(t))$, $0\le t<T_1(\bar t_2)$, has the same
distribution as a chordal SLE$(\kappa_1;-\frac{\kappa_1}2,\vec{\rho}_1)$ trace started
from $(\vphi_2(\bar t_2,x_1);\xi_2(\bar t_2), \vphi_2(\bar t_2,\vec{p}))$. Then we find
that a.s.\ $\lim_{t\to T_1(\bar t_2)} \vphi_2(\bar t_2,\gamma_1(t))\in (\xi_2(\bar t_2),\vphi_2(\bar t_2,p_2))$.
Thus $\vphi_2(\bar t_2,\gamma_1(t))$, $0\le t<T_1(\bar t_2)$, disconnects $\xi_2(\bar t_2)$ from
$\infty$ in $\HH$, and so $\gamma_1$ disconnects $\gamma_2(\bar t_2)$ from $\infty$ in $\HH\sem L_2(\bar t_2)$. By choosing a sequence of
$(\F^2_t)$-stopping times that are dense in $[0,T_2)$, we conclude that a.s.\
$\lin{K_2(T_2^-)}\subset\lin{D_\HH(\gamma_1)}$,
where $K_2(T_2^-)=\cup_{0\le t<T_2} K_2(t)$.  From
Proposition \ref{more force} and Proposition \ref{coordinate}, a.s.\ $x_1$ is a subsequential limit of $\gamma_2(t)$ as $t\to T_2$.
Similarly, for every $(\F^1_t)$-stopping time $\bar t_1\in(0,\bar T_1)$, $\gamma_1(\bar t_1)$ is
a subsequential limit of $\gamma_2(t)$ as $t\to T_2(\bar t_1)$. By choosing a sequence of
$(\F^1_t)$-stopping times that are dense in $[0,T_1)$, we conclude that a.s. $\gamma_1(t)\in
 \lin{K_2(T_2^-)}$ for $0\le t< T_1$. So we have the following lemma and theorem. Here  $\pa^{\out}_\HH S$ is defined for
bounded $S\subset\HH$, which is the intersection of $\HH$ with the boundary of the
unbounded component of $\HH\sem S$. For detailed proof of the lemma, please see Lemma 5.1 in \cite{duality}.

\begin{Lemma} Almost surely
$\pa^{\out}_{\HH}K_2(T^-_2)$ is the image of $\gamma_1(t)$, $0<t<T_1$. \label{K-infty}
\end{Lemma}

\begin{Theorem} Suppose $\kappa>4$;
 $p_1<x_1<p_3<p_4<x_2<p_2$; $C_2,C_4\ge 1/2$, $C_1=1-C_2$, and $C_3=1-C_4$. Let
 $K(t)$, $0\le t<T$, be chordal SLE$(\kappa;-\frac{\kappa}2,C_1(\kappa-4),C_2(\kappa-4),C_3(\kappa-4),
C_4(\kappa-4))$ process
 started from $(x_2;x_1,p_1,p_2,p_3,p_4)$. Let $K(T^-)=\cup_{0\le t<T} K(t)$.
 Then a.s.\ $K(T^-)$ is bounded, and $\pa^{\out}_\HH K(T^-)$
  has the distribution of the image of a chordal
 SLE$(\kappa';-\frac{\kappa'}2,C_1(\kappa'-4),C_2(\kappa'-4),C_3(\kappa'-4),C_4(\kappa'-4))$
 trace started from $(x_1;x_2,p_1,p_2,p_3,p_4)$, where
 $\kappa'=16/\kappa$. \label{K-infty*}
\end{Theorem}

The above lemma and theorem still hold if we let $p_1\in(-\infty,x_1)$,
or $=x_1^-$; let $p_2\in(x_2,\infty)$, or $=\infty$,
or $=x_2^+$; let $p_3\in(x_1,x_2)$, or $=x_1^+$; let $p_4\in(x_1,x_2)$ or $=x_2^-$. Here
if $p_2=x_2^+$, we use Theorem \ref{tendstoinfty} instead of Theorem \ref{lim-1p3}
to prove that the image of $\gamma_1$ in Lemma \ref{K-infty} is a crosscut in $\HH$ on $\R$.

\vskip 3mm

\no {\bf Proof of Theorem \ref{kappa(4,8)swallow}.} First suppose $x<0$. Then $\lambda$ is supported by $(-\infty,x)$, and
$\bar\mu=\int \bar\nu_y\,d\lambda(y)$ follows from Corollary \ref{kappa(4,8)} and Theorem
\ref{K-infty*} with $x_1=y$, $x_2=0$, $p_1=y^-$, $p_2=\infty$, $p_3=y^+$, $p_4=x$,
$C_1=\frac{\kappa-6}{\kappa-4}$, $C_2=\frac{2}{\kappa-4}$,
$C_3=-1/2$, and $C_4=1$. From Theorem \ref{lim-1p3} and Proposition \ref{coordinate2},
for each $y\in(-\infty,x)$, $\bar\nu_y$ is supported by the space of crosscuts in $\HH$ from $y$ to
some point on $(0,\infty)$.
Thus a.s.\
$\pa^{\out}_\HH K(T_x)$ is a crosscut in $\HH$ on $\R$ connecting some $y\in(-\infty,x)$ with some $z\in(0,\infty)$.
The case that $x>0$ is symmetric. $\Box$

\vskip 3mm

\vskip 3mm Let $S\subset\HH$. Suppose
$\lin{S}\cap(a,\infty)=\emptyset$ for some $a\in\R$. Then there is a
unique component of $\HH\sem\lin{S}$, which has $(a,\infty)$ as part
of its boundary. Let $D_+$ denote this component. Then $\pa D_+\cap
\HH$ is called the right boundary of $S$ in $\HH$. Let it be denoted
by $\pa_\HH^+ S$. Similarly, if $\lin{S}\cap(-\infty,a)=\emptyset$
for some $a\in\R$. Then there is a unique component of $\HH\sem\lin{S}$,
which has $(-\infty,a)$ as part
of its boundary. Let $D_-$ denote this component. Then $\pa D_-\cap
\HH$ is called the left boundary of $S$ in $\HH$. Let it be denoted
by $\pa_\HH^- S$. The following theorem improves Theorem 5.3 in \cite{duality}.

\begin{Theorem} Let $\kappa>4$ and $C_r,C_l\ge 1/2$. Let
$K(t)$, $0\le t<\infty$, be a
chordal SLE$(\kappa;C_r(\kappa-4),C_l(\kappa-4))$ process started
from $(0;0^+,0^-)$. Let $K(\infty)=\cup_{t\ge 0}K(t)$. Let
 $\kappa'=16/\kappa$ and $W(z)=1/\lin{z}$. Then (i) $W(\pa_\HH^+ K(\infty))$ has the same distribution
as the image of a chordal SLE$(\kappa';(1-C_r)(\kappa'-4),(1/2-C_l)(\kappa'-4))$ trace started
from $(0;0^+,0^-)$; (ii) $W(\pa_\HH^- K(\infty))$ has the same distribution
as the image of a chordal SLE$(\kappa';(1/2-C_r)(\kappa'-4),(1-C_l)(\kappa'-4))$ trace started
from $(0;0^+,0^-)$; and (iii) a.s.\ $\lin{K(\infty)}\cap\R=\{0\}$.
\label{duality1}
\end{Theorem}
{\bf Proof.}
Let $W_0(z)=1/(1-z)$. Then
$W_0$ maps $\HH$ conformally onto $\HH$, and $W_0(0)=1$, $W_0(\infty)=0$, $W_0(0^\pm)=1^\pm$.
From Proposition \ref{coordinate}, after a time-change, $(W_0(K(t)))$ has the same
distribution as a chordal SLE$(\kappa;(\frac 32-C_r-C_l)(\kappa-4)-\frac\kappa
2,C_r(\kappa-4),C_l(\kappa-4))$ process started from
$(1;0,1^+,1^-)$. Applying Theorem \ref{K-infty*} with $x_1=0$,
$x_2=1$, $p_1=0^-$, $p_2=1^+$, $p_3=0^+$, $p_4=1^-$,
$C_1=1-C_r$, $C_2=C_r$, $C_3=1/2-C_l$, and $C_4=C_l$, we
find that $\pa_\HH^{\out} W_0(K(\infty))$ has the same distribution
as the image of a chordal SLE$(\kappa';(C_2+C_4)(\kappa'-4)-\frac{\kappa'} 2,
C_1(\kappa'-4),C_3(\kappa'-4))$ trace started from $(0;1,0^-,0^+)$. Let
$\gamma$ denote this trace. From Proposition \ref{coordinate} and Theorem
\ref{tendstoinfty}, $\gamma$ is a crosscut in $\HH$ from $0$ to $1$. Thus
$\pa_{\HH}^{+} K(\infty)= W_0^{-1}(\gamma)$, and so $W(\pa_{\HH}^{+}
K(\infty))=W\circ W_0^{-1}(\gamma)$. Let $W_1=W\circ W_0^{-1}$. Then
$W_1(z)=\lin{z}/(\lin{z}-1)$. So $W_1(0)=0$,
$W_1(1)=\infty$, $W_1(0^\pm)=0^\mp$. From Proposition \ref{coordinate}, after a time-change,
$W_1(\gamma)$ has the same distribution as  a chordal
SLE$(\kappa';C_1(\kappa'-4),C_3(\kappa'-4))$ trace started from $(0;0^+,0^-)$.
Since $C_1=1-C_r$ and $C_3=1/2-C_l$, so we have (i). Now (ii) follows from symmetry.
Finally, from (i), (ii), and Proposition \ref{3pt**}, $\pa_{\HH}^{+} K(\infty)$ and $\pa_{\HH}^{-} K(\infty)$
are two crosscuts in $\HH$ that connect $\infty$ with $0$, so we have (iii).
$\Box$

\vskip 3mm

In the proof of the above theorem, if we choose $p_2$ and $p_4$ to be generic force points,
then we may obtain the following theorem using a similar argument.

\begin{Theorem} Let $\kappa>4$, $C_r,C_l\ge 1/2$, and $p_r>0>p_l$.
Suppose $K(t)$, $0\le t<\infty$, is a chordal SLE$(\kappa;C_r(\kappa-4),C_l(\kappa-4))$
process started from $(0;p_r,p_l)$. Let $K(\infty)=\cup_{t\ge 0} K(t)$ and $\kappa'=16/\kappa$.
Then $\pa_\HH^+ K(\infty)$ is a
crosscut in $\HH$ from $\infty$ to some point on $(0,p_r)$; $\pa_\HH^- K(\infty)$
is a crosscut in $\HH$ from $\infty$ to some point on $(p_l,0)$; and $K(\infty)$ is bounded away
from $(-\infty,p_l]$ and $[p_r,+\infty)$. \label{leftrightboundary}
\end{Theorem}

\section{Boundary of Chordal SLE} \label{BoundaryofChordalSLE}
In this section, we use Theorem \ref{kappa(4,8)swallow} and Proposition \ref{kappa>8}
to study the boundary of standard chordal SLE$(\kappa)$ hulls for $\kappa>4$.

Let $\kappa>4$. Let $K(t)$, $0\le t<\infty$, be a standard chordal SLE$(\kappa)$ process. Let $\xi(t)$,
 $\vphi(t,\cdot)$, and $\gamma(t)$, $0\le t<\infty$, be the corresponding
driving function, chordal Loewner maps, and trace. Then there is a Brownian motion $B(t)$ such that
$\xi(t)=\sqrt\kappa B(t)$, $t\ge 0$.
For each $t>0$, let $a(t)=\inf (\lin{K(t)}\cap\R)$ and $b(t)
=\sup (\lin{K(t)}\cap\R)$, then $a(t)<0<b(t)$, and $\vphi(t,\cdot)$ maps $(-\infty,a(t))$ and
$(b(t),+\infty)$ onto $(-\infty,c(t))$ and $(d(t),+\infty)$ for some $c(t)<0<d(t)$. And we have
$c(t)\le\xi(t)\le d(t)$, $t>0$. For each $t>0$, $f_t:=\vphi(t,\cdot)^{-1}$ extends continuously to
$\lin{\HH}$ with $f_t(c(t))=a(t)$, $f_t(d(t))=b(t)$, $f_t(\xi(t))=\gamma(t)$, and
 $ K(t)$ is bounded by $f_t([c(t),d(t)])$ and $\R$. We have the following theorem.

\begin{Theorem}
Let $T\in(0,\infty)$ be a stopping time w.r.t.\ the filtration generated by $(\xi(t))$.
Then $\gamma(T)\in\R$ a.s.\ implies that $\xi(T)=c(T)$ or $=d(T)$, and the curve
$f_t(x)$, $c(T)< x< d(T)$, is a crosscut in $\HH$ on $\R$ with dimension $1+2/\kappa$ everywhere; and
$\gamma(T)\in\HH$ a.s.\ implies that $c(T)<\xi(T)<d(T)$, and the two curves
$f_t(x)$, $c(T)< x\le \xi(T)$, and  $f_t(x)$, $\xi(T)\le x< d(T)$, are both semi-crosscuts in $\HH$
on $\R$ with dimension $1+2/\kappa$ everywhere. Moreover, $\lin{K(T)}$ is connected, and has no cut-point on $\R$. \label{boundaryofSLE}
\end{Theorem}

Here a curve $\alpha$ is said to have dimension $d$ everywhere if any non-degenerate subcurve of $\alpha$
has Hausdorff dimension $d$. From the main theorem in \cite{dim-SLE}, every standard chordal
SLE$(\kappa)$ trace has dimension $(1+\kappa/8)\wedge 2$ everywhere. From Girsanov Theorem and Proposition
\ref{coordinate2},
this is also true for any chordal or strip SLE$(\kappa;\vec{\rho})$ trace. For a connected set $K\subset\C$, $z_0\in K$
is called a cut-point of $K$, if $K\sem \{z_0\}$ is not connected. Such cut-point must lie on the boundary of $K$.

We need a lemma to prove this theorem. For each $p\in\R\sem\{0\}$, let $T_p$ denote the first
time that $p$ is swallowed by $K(t)$. Then $T_p>0$ is a finite stopping time because
$\kappa>4$.

\begin{Lemma}
For $p_-<0<p_+$, the events $\{T_{p_-}<T_{p_+}\}$ and $\{T_{p_+}<T_{p_-}\}$ both have positive
probabilities. \label{compare}
\end{Lemma}
{\bf Proof.} Let $T=T_{p_-}\wedge T_{p_+}$. Let $X_\pm(t)=\vphi(t,p_\pm)-\xi(t)$, $0\le t<T$.
Then $X_\pm(t)$ satisfies the SDE: $dX_\pm(t)=-\sqrt\kappa dB(t)+\frac{2}{X_\pm(t)}dt$.
Let $Y_\pm(t)=\ln(|X_\pm(t)|)$, $0\le t<T$. From Ito's formula,  $Y_\pm(t)$ satisfies
the SDE:
$$dY_\pm(t)=-\frac{\sqrt\kappa}{X_\pm(t)} dB(t)+\Big(2-\frac\kappa 2\Big)\,\frac{dt}{X_\pm(t)^2}.$$
Let $Y(t)=Y_+(t)-Y_-(t)$, $0\le t<T$. Then $Y(t)$ satisfies the SDE:
$$dY(t)=-\sqrt\kappa\Big[\frac 1{X_+(t)}-\frac 1{X_-(t)}\Big]\,dB(t)+\Big(2-\frac\kappa 2\Big)\,
\Big[\frac{1}{X_+(t)^2}-\frac{1}{X_-(t)^2}\Big]\,dt.$$
Let $u(t)=\int_0^t (1/X_+(s)-1/X_-(s))^2ds$, $0\le t<T$. Let $Z(t)=Y(u^{-1}(t))$, $0\le t<u(T)$.
Then there is a Brownian motion $\til B(t)$ such that $Z(t)$ satisfies the SDE:
$$dZ(t)=-\sqrt\kappa d\til B(t)+\Big(2-\frac\kappa 2\Big)\,\frac{X_-(u^{-1}(t))+X_+(u^{-1}(t))}
{X_-(u^{-1}(t))-X_+(u^{-1}(t))}\,dt$$$$=-\sqrt\kappa d\til B(t)+\Big(\frac\kappa 2-2\Big)\,
\tanh_2(Z(t))\,dt.$$
From the chordal Loewner equation, $X_+(t)-X_-(t)=\vphi(t,p_+)-\vphi(t,p_-)$ increases in $t$. If $T=T_{p_-}$,
as $t\to T^-$, $X_-(t)=\vphi(t,p_-)-\xi(t)\to 0$, so $|X_+(t)|/|X_-(t)|\to \infty$, which implies that
$Z(t)\to +\infty$ as $t\to u(T)$. Similarly, if $T=T_{p_+}$, then $Z(t)\to -\infty$ as $t\to u(T)$.
Thus as $t\to T$, either $Z(t)\to +\infty$ or $Z(t)\to -\infty$.
For $x\in\R$, let $h(x)=\int_0^x \cosh_2(s)^{2/\kappa-2}ds$. Since $2/\kappa-2<0$, so $h$ maps $\R$ onto a
finite interval, say $(-L,L)$. And we have $\frac\kappa 2 h''(x)+(\frac\kappa 2-2)h'(x)\tanh_2(x)=0$ for any $x\in\R$.
Let $W(t)=h(Z(t))$, $0\le t<u(T)$. Then as $t\to u(T)$, either $W(t)\to L$ or $W(t)\to -L$.
From Ito's formula, $(W(t))$ is a bounded martingale. Thus the probability that $\lim_{t\to u(T)}W(t)=L$
is $(W(0)-(-L))/(2L)>0$. So the probability that $T_{p_-}<T_{p_+}$, i.e., $T=T_{p_-}$, is positive.
Similarly, the probability that $T_{p_+}<T_{p_-}$  is also positive. $\Box$

\vskip 3mm

\no{\bf Proof of Theorem \ref{boundaryofSLE}.} Let $\kappa'=16/\kappa\in(0,4)$.
 If $T=T_p$ for some $p\in\R\sem\{0\}$, then $\gamma(T)\in\R$, and
 $\xi(T)= c(T)$ or $d(T)$, depending on
whether $p<0$ or $p>0$. From Theorem \ref{kappa(4,8)swallow} and Proposition \ref{kappa>8},
$\pa  K(T)\cap\HH=\{f_{T}(x):c(T)<x<d(T)\}$ is the image of a
chordal SLE$(\kappa',\vec{\rho})$ trace, and so it
has dimension $1+\kappa'/8=1+2/\kappa$ everywhere. We also see that
this curve is a crosscut in $\HH$ on $\R$, so $K(T)$ is the hull bounded by this crosscut. Thus
$\lin{K(T)}$ is connected, and has no cut-point.

Now consider the general case.
We first prove (i):  $\xi(T)=c(T)$ a.s.\ implies that
$f_t(x)$, $c(T)< x< d(T)$, is a crosscut in $\HH$ on $\R$ with dimension $1+2/\kappa$ everywhere.
Let ${\cal E}$ denote the event that  $\xi(T)=c(T)$, but $f_t(x)$, $c(T)< x < d(T)$, is not a crosscut
in $\HH$ on $\R$, or does not have dimension $1+2/\kappa$ everywhere. Assume that $\PP({\cal E})>0$.
For each $n\in\N$, let
$${\cal E}_n:= \{\xi(T)=c(T)\}\cap\{-n<a(T)\}\cap\{d(T)-c(T)>1/n\}\cap$$
$$\cap \{f_t(x),c(T)+1/n\le x< d(T) \mbox{, is
not a semi-crosscut in }\HH\mbox{ on }\R,$$
$$\mbox{or does not have dimension }1+2/\kappa\mbox{ everywhere}\}.$$
Since $f_T(c(T))=a(T)\in\R$, and $a(T)< b(T)=f_T(d(T))$, so ${\cal E}=\cup_{n=1}^\infty{\cal E}_n$. Then
 there is $n_0\in\N$ such that $\PP({\cal E}_{n_0})>0$.

Let $(\til K(t),0\le t<\infty)$ be a standard chordal SLE$(\kappa)$ process that is independent of $(K(t))$.
Let $\til{\cal E}_{n_0}$ denote the event that $\til K(t)$ swallows $\vphi(T,-n_0)-\xi(T)$ before swallowing $1/n_0$,
and let $\til{T}$ denote the first time that $\til K(t)$ swallows $\vphi(T,-n_0)-\xi(T)$.
From Lemma \ref{compare}, the probability of $\til{\cal E}_{n_0}$ is positive.
Let $\ha{\cal E}_{n_0}={\cal E}_{n_0}\cap \til{\cal E}_{n_0}$. Then $\ha{\cal E}_{n_0}$ also has
positive probability.

Define $\ha K(t)=K(t)$ for $0\le t\le T$; and $\ha K(t)=K(T)\cup f_T(\til K(t-T)+\xi(T))$ for $t>T$. Then
$(\ha K(t))$ has the same distribution as $(K(t))$. Let $\ha{T}_{-n_0}$ denote the first time
that $\ha K(t)$ swallows $-n_0$. Then  $\pa \ha K(\ha T_{-n_0})\cap\HH$ is a.s.\ a crosscut in $\HH$
on $\R$ with dimension $1+2/\kappa$ everywhere. Since on $\ha{\cal E}_{n_0}$, $\ha{T}_{-n_0}
=T+\til T$, and  $\lin{\til K(\til T)}\cap\R$ is bounded above by $1/n_0$, so
$\{f_T(x), c(T)+1/n_0\le x< d(T)\}$ is a subset of the boundary of
$\ha K(\ha T_{-n_0})=K(T)\cup f_T(\til K(\til T)+\xi(T))$ in $\HH$,
which implies that a.s.\  $f_T(x)$, $c(T)+1/n_0\le x< d(T)$, is a semi-crosscut
with dimension $1+2/\kappa$ everywhere. This contradicts that $\ha{\cal E}_{n_0}$  has positive
probability. So we have (i). Symmetrically, we have (ii):
$\xi(T)=d(T)$ a.s.\ implies that
$f_t(x)$, $c(T)< x< d(T)$, is a crosscut in $\HH$ on $\R$ with dimension $1+2/\kappa$ everywhere.

If $\gamma(T)=f_T(\xi(T))\in\HH$, then $\gamma(T)\not\in\{c(T),d(T)\}$, so $c(T)<\xi(T)<d(T)$.
Using the same argument as in (i), we can prove (iii): $\gamma(T)\in\HH$ a.s.\ implies
that $f_t(x)$, $\xi(T)\le x < d(T)$, is a semi-crosscut in $\HH$ on $\R$ with dimension $1+2/\kappa$ everywhere.
Symmetrically, we have (iv): $\gamma(T)\in\HH$ a.s.\ implies
that $f_t(x)$, $c(T)< x\le \xi(T)$, is a semi-crosscut in $\HH$ on $\R$ with dimension $1+2/\kappa$ everywhere.
From (iii) and (iv), we see that $\gamma(T)\in\HH$ a.s.\ implies that $\lin{K(T)}$ is connected, and has no cut-point
on $\R$.  Similarly, we have (v): $c(T)<\xi(T)<d(T)$ and $\gamma(T)\in\R$ a.s.\ implies that
$f_T(x)$, $\xi(T)< x< d(T)$, and $f_T(x)$, $c(T)< x < \xi(T)$, are both crosscuts or degenerate
crosscuts in $\HH$ on $\R$. Moreover, these two curves intersect
at only one point: $\gamma(T)$, since the curve $\alpha(y):=f_T(\xi(T)+iy)$, $y> 0$,
connects $\gamma(T)$ with $\infty$, and does not intersect the above two curves. So $\gamma(T)$ is
a cut-point of $\lin{K(T)}$ on $\R$.

To finish the proof, it remains to prove (vi): $\gamma(T)\in\R$ a.s.\ implies that $\xi(T)=c(T)$ or $=d(T)$.
Let ${\cal E}$ denote the event that $\gamma(T)\in\R$ and $c(T)<\xi(T)<d(T)$.
We suffice to show that $\PP({\cal E})=0$. Assume that $\PP({\cal E})>0$. Assume that $\cal E$ occurs.
From (v), we know that
$K(T)=K_1\cup K_2$, where $K_1$ and $K_2$ are hulls bounded by crosscut or degenerate crosscut in $\HH$ on $\R$,
and $\lin{K_1}\cap \lin{K_2}=\{\gamma(T)\}$. Since $\kappa>4$, so a.s.\ $K(T)$ contains a neighborhood of
$0$ in $\HH$. We may label $K_1$ and $K_2$ such that $K_1$ contains a neighborhood of $0$ in $\HH$. Then
$\gamma(T)\ne 0$.
Let ${\cal S}=\{\lin{\B(x+iy;r)}:x,y,r\in\Q, y,r>0, r<y/2\}$, where $\B(z_0;r):=\{z\in\C:|z-z_0|<r\}$.
Then $\cal S$ is countable, and every $A\in{\cal S}$ is contained in $\HH$. For $A\in{\cal S}$,
let ${\cal E}_A$ denote the intersection of $\cal E$ with the event that $A\cap \pa  K_2
\not=\emptyset$ and $A\cap K_1=\emptyset$. Then ${\cal E}=\cup_{A\in{\cal S}}{\cal E}_A$. So
there is $A_0\in\cal S$ such that $\PP({\cal E}_{A_0})>0$. Let $T_0$ be the first time that $\gamma(t)$
hits ${A_0}$. Let $T_1=T\wedge T_0$. Then $T_1$ is a finite stopping time. Assume ${\cal E}_{A_0}$ occurs.
Since $\gamma(t)$, $0\le t\le T$,
visits every point on $\pa  K_2\cap\HH\subset \pa K(T)\cap\HH$, so
$T_0\le T$, and so $T_1=T_0$. We have $\gamma(T_1)=\gamma(T_0)\in {A_0}$.
Since ${A_0}\cap\R=\emptyset$, so $\gamma(T_1)\in\HH$.
Since $\gamma(0)=0\in\lin{K_1}$, and $\gamma(T_1)\in \lin{K_2}$, which are both different from $\gamma(T)$,
so $\gamma(T)\in\lin{K_1}\cap \lin{K_2} $ is a cut-point of $\lin{K(T_1)}$. However, since $T_1$ is a positive
finite stopping time, and $\gamma(T_1)\in\HH$ on ${\cal E}_{A_0}$, so
from (iii) and (iv) in the above proof, a.s.\ $\lin{K(T_1)}$ has no cut-point on $\R$  in the event ${\cal E}_{A_0}$. This contradicts that
$\PP({\cal E}_{A_0})>0$. So $\PP({\cal E})=0$. $\Box$

\begin{Corollary} For any stopping time $T\in(0,\infty)$, a.s.\ $f_T(x)\not\in\R$ for $x\in(c(T),d(T))$;
$\pa  K(T)\cap\HH$ has Hausdorff dimension $1+2/\kappa$; $\lin{K(T)}$ is connected, and has no cut-point on $\R$;
  and for every $x\in(a(T),b(T))$, $K(T)$ contains a neighborhood of $x$ in $\HH$.
\end{Corollary}

In the above theorem, when $\gamma(T)\in\HH$, $\pa K(T)\cap\HH$ is composed of two semi-crosscuts in $\HH$ on $\R$, which are
$f_T(x)$, $c(T)< x\le \xi(T)$, and $f_T(x)$, $\xi(T)\le x< d(T)$. If the two semi-crosscuts intersect
only at $\gamma(T)=f_T(\xi(T))$, then we get a crosscut $f_T(x)$, $c(T)< x < d(T)$. If the
two semi-crosscuts intersect at any point $z_0$ other than $\gamma(T)$, then $z_0$ is a cut-point
of $K(T)$. To see this, suppose $f_T(x_1)=f_T(x_2)=z_0$, where $c(T)<x_1<\xi(T)<x_2<d(T)$. Then
$f_T(x)$, $c(T)<x\le x_1$, and $f_T(x)$, $x_2\le x< d(T)$, are two semi-crosscuts in $\HH$ on $\R$,
which together bound a hull in $\HH$ on $\R$. Let it be denoted by $K_1$.
The simple curves
$f_T(x)$, $x_1\le x\le \xi(T)$, and $f_T(x)$, $\xi(T)\le x\le x_2$, together bound a closed bounded set
in $\HH$. Let it be denoted by $K_2$. Then $K(T)=K_1\cup K_2$ and $K_1\cap K_2=\{z_0\}$. On the other
hand, every cut-point of $K(T)$ corresponds to an intersection point between
$f_T(x)$, $c(T)< x< \xi(T)$, and $f_T(x)$, $\xi(T)< x< d(T)$, and so such cut-point disconnects
$\gamma(T)$ from $\xi(0)=0$ in $K(T)$. From Theorem 5 in \cite{cut-SLE},
if $\kappa>8$ and $T>0$ is a constant, then a.s.\ $K(T)$ has no cut-point, so
$f_T(x)$, $c(T)< x< d(T)$, is a crosscut in $\HH$ on $\R$. We now make some improvement over this
result.

\begin{Theorem}
If $\kappa\ge 8$ and $T\in(0,\infty)$ is a stopping time, then a.s.\ $K(T)$ has no cut-point, and so
$f_T(x)$, $c(T)< x< d(T)$, is a crosscut in $\HH$ on $\R$.
\end{Theorem}
{\bf Proof.} First suppose $\kappa>8$.
Let $\cal E$ denote the event that $K(T)$ has a cut-point. We suffice
to show that $\PP({\cal E})=0$. Assume that $\PP({\cal E})>0$. For each $n\in\N$, let
${\cal E}_n$ denote the event that $c(T)+1/n<\xi(T)<d(T)-1/n$, and the two curves
$f_T(x)$, $c(T)< x \le\xi(T)-1/n$, and $f_T(x)$, $\xi(T)+1/n\le x< d(T)$, are not disjoint.
Then ${\cal E}=\cup_{n=1}^\infty {\cal E}_n$. So there is $n_0\in\N$ such that $\PP({\cal E}_{n_0})>0$.

Let $(\til K(t))$ be a standard chordal SLE$(\kappa)$ process that is independent of $(K(t))$.
There is a small $h>0$ such that the probability that $\lin{\til K(h)}\cap\R\subset(-1/n_0,1/n_0)$ is positive.
There is $t_0\in[0,\infty)$ such that $\PP({\cal E}_{n_0}\cap \{t_0-h\le T\le t_0\})>0$.
Let $\ha {\cal E}$ denote the intersection of ${\cal E}_{n_0}\cap \{t_0-h\le T\le t_0\}$ with
$\{\lin{\til K(h)}\cap\R\subset(-1/n_0,1/n_0)\}$. Then $\ha {\cal E}$ also has positive probability.
Define $\ha K(t)=K(t)$ for $0\le t\le T$; and $\ha K(t)=K(T)\cup f_T(\til K(t-T)+\xi(T))$ for $t>T$. Then
$(\ha K(t))$ has the same distribution as $(K(t))$. From Theorem 5 in \cite{cut-SLE}, a.s.\ $\ha K(t_0)$ has no cut-point.
Since $T\le t_0\le T+h$, so $K(T)\subset \ha K(t_0)\subset K(T)\cup f_T(\til K(h)+\xi(T))$.
In the event $\ha {\cal E}$, since $\lin{\til K(h)}\cap\R\subset(-1/n_0,1/n_0)$, so
$f_T(x)$, $c(T)< x \le\xi(T)-1/n_0$, and $f_T(x)$, $\xi(T)+1/n_0\le x< d(T)$, are subarcs of
$\pa  \ha K(t_0)\cap\HH$. However, in the event $\ha {\cal E}$,  the above two curves are not disjoint,
so $\ha K(t_0)$ has a cut-point, which contradicts that $\ha {\cal E}$ has positive probability. Thus $\PP({\cal E})=0$.

Now suppose $\kappa=8$. Let $\gamma^R(t)=\gamma(1/t)$, $0<t<\infty$. Since chordal SLE$(8)$ trace
is reversible (c.f.\ \cite{LSW-2}), so after a time-change, $\gamma^R$ has the distribution of a
chordal SLE$(8)$ trace in $\HH$ from $\infty$ to $0$. Thus a.s.\ there is a crosscut $\alpha$
in $\HH\sem \gamma^R((0,1/T])=\HH\sem\gamma([T,\infty)$ connecting $\gamma^R(1/T)=\gamma(T)$ with $0$.
Then $\alpha\subset K(T)$ and does not intersect $\pa K(T)$. If $K(T)$ has any cut-point, the cut-point
must disconnect $\gamma(T)$ from $0$ in $K(T)$, so such $\alpha$ does not exist. Thus a.s.\ $K(T)$ has no
cut-point. $\Box$

\vskip 3mm

If $\kappa\in(4,8)$, this theorem does not hold since from Theorem 5 in \cite{cut-SLE}, the probability
that $K(1)$ has cut-point is positive.

\section{More Geometric Results} \label{MoreGeometricResults}
The description of the boundary of SLE$(\kappa)$ hulls for $\kappa>4$ enables us to obtain some
results about the limit of SLE$(\kappa;\vec{\rho})$ traces when $\kappa>4$. We will prove that the limits of the traces exist
when certain conditions are satisfied.

Let $\kappa>4$. In this section, $L(t)$, $0\le t<T_e$, is a strip SLE$(\kappa;\vec{\rho})$ process started
from $(0;\vec{p})$, where no force point is degenerate. Let $\xi(t)$, $\psi(t,\cdot)$, and $\beta(t)$,
$0\le t<T_e$, be the corresponding driving function, strip Loewner maps, and trace.
For $t\in(0,T_e)$, let $a(t)=\inf(\lin{L(t)}\cap\R)<0$ and $b(t)=\sup(\lin{L(t)}\cap\R)>0$.
Then $\psi(t,\cdot)$ maps $(-\infty,a(t))$ and
$(b(t),+\infty)$ onto $(-\infty,c(t))$ and $(d(t),+\infty)$ for some $c(t)<0<d(t)$, and we have
$c(t)\le\xi(t)\le d(t)$. For each $t>0$, $f_t:=\psi(t,\cdot)^{-1}$ extends continuously to
$\lin{\St_\pi}$ such that $f_t(c(t))=a(t)$, $f_t(d(t))=b(t)$, and $f_t(\xi(t))=\beta(t)$. From
Theorem \ref{boundaryofSLE}, Proposition \ref{coordinate2}, and Girsanov Theorem, we have the
following lemma.

\begin{Lemma}
If $T\in(0,T_e)$ is a stopping time, then a.s.\ $f_T(x)\in\St_\pi$ for $c(T)<x<d(T)$, and for
every $x\in(a(T),b(T))$, $L(T)$ contains a neighborhood of $x$ in $\St_\pi$. \label{boundaryofstrip}
\end{Lemma}

\begin{Lemma} Let $T\in[0,T_e)$ be a stopping time. Define $\beta_T(t)=\psi(T,\beta(T+t))-\xi(T)$,
$0\le t<T_e-T$. Suppose $\vec{p}=(p_1,\dots,p_N)$. If $\psi(T,p_m)-\xi(T)=p_m$ for $1\le m\le N$,
then $\beta_T$ has the same distribution as $\beta$. In the general case, conditioned on $\beta(t)$,
$0\le t\le T$, $\beta_T$ is a strip SLE$(\kappa;\vec{\rho})$ trace started from $(0;\vec{q})$,
where $\vec{q}=(q_1,\dots,q_N)$ and $q_m=\psi(T,p_m)-\xi(T)$, $1\le m\le N$.
\label{betaT}
\end{Lemma}
{\bf Proof.} This follows from the definition of strip SLE$(\kappa;\vec{\rho})$ process and the
property that Brownian motion has i.i.d.\ increment. $\Box$

\begin{Lemma}
Let $\kappa>4$, $\rho_+,\rho_-\in\R$, $\rho_++\rho_-=\kappa-6$, and $\rho_--\rho_+\ge 2$.
Suppose $\beta(t)$, $0\le t<\infty$, is a strip SLE$(\kappa;\rho_+,\rho_-)$ trace started from
$(0;+\infty,-\infty)$. Then a.s.\ any subsequential limit of $\beta(t)$ as $t\to \infty$ does not
lie on $\R\cup\R_\pi\cup\{-\infty\}$. \label{nootherlimit}
\end{Lemma}
{\bf Proof.} Let $Q$ denote the set of subsequential limits of $\beta(t)$ as $t\to\infty$.
 Let $\sigma=(\rho_--\rho_+)/2\ge 1$. Then there is a Brownian motion $B(t)$ such that
$\xi(t)=\sqrt\kappa B(t)+\sigma t$, $0\le t<\infty$. Thus a.s.\ there is a random number $A_0<0$ such that
 $\xi(t)\ge A_0$ for $0\le t<\infty$. From (\ref{strip}), for any $z\in\St_\pi$ with
$\Ree z<A_0$,   $\psi(t,z)$ never blows up for $0\le t<\infty$. Thus a.s.\
$\beta([0,\infty))\subset \{z\in\lin{\St_\pi}:\Ree z\ge A_0\}$. So a.s.\ $-\infty\not\in Q$.
Moreover, for any $\eps>0$, there is $R_\eps>0$ such that the probability that $\Ree\beta(t)\ge -R_\eps$ for
$0\le t<\infty$ is at least $1-\eps$.

Fix $x_0\in\R$. Let $X(t)=\Ree \psi(t,x_0+\pi i)-\xi(t)$, $0\le t<\infty$. Then $X(t)$ satisfies
the SDE: $dX(t)=-\sqrt\kappa dB(t)+\tanh_2(X(t))dt-\sigma dt$. Define $h$ on $\R$ such that
$$h'(x)=\exp(x/2)^{\frac 4\kappa\,\sigma}\cosh_2(x)^{-\frac 4\kappa},\quad x\in\R.$$
Since $\sigma\ge 1$, so $h$ maps $\R$ onto $(L,\infty)$ for some $L\in\R$.
Let $Y(t)=h(X(t))$, $0\le t<\infty$. From Ito's formula, $Y(t)$ satisfies the SDE:
$dY(t)=-h'(X(t))\sqrt\kappa dB(t)$. Define $u(t)=\int_0^t \kappa h'(X(s))^2 ds$, $0\le t<\infty$,
and $u(\infty)=\sup u([0,\infty))$.
Then $Y(u^{-1}(t))$, $0\le t<u(\infty)$, has the distribution of a partial Brownian motion.
Since $Y(u^{-1}(t))\in(L,\infty)$ for $0\le t<u(\infty)$, so a.s.\ $u(\infty)<\infty$ and
$\lim_{t\to\infty} Y(t)=\lim_{t\to u(\infty)} Y(u^{-1}(t))\in [L,\infty)$. Note that
$\lim_{t\to\infty} Y(t)\in (L,\infty)$ implies that $\lim_{t\to\infty}X(t)\in\R$ and so
$X(t)$, $0\le t< \infty$, is bounded. If $X$ is bounded on $[0,\infty)$,
from the definition of $u$, $u'(t)$ is uniformly bounded below by a positive constant,
which implies that $u(\infty)=\infty$.
Since a.s.\ $u(\infty)<\infty$, so $\lim_{t\to\infty} Y(t)\not\in (L,\infty)$.
Thus a.s.\ $\lim_{t\to\infty} Y(t)=L$, and so $\lim_{t\to\infty} X(t)=-\infty$.

Fix $\eps>0$. Let $T$ be the first time such that $X(t)\le -R_\eps-1$.
Then $T$ is a finite stopping time. Let $\beta_T$ be defined as in Lemma \ref{betaT}.
Then $\beta_T$ has the same distribution as $\beta$.
So the probability that $\Ree \beta_T(t)\ge -R_\eps$ for
any $0\le t<\infty$ is at least $1-\eps$. Let $Q_T$ denote the set of subsequential limits
of $\beta_T(t)$ as $t\to\infty$. Then the probability that $Q_T\cap(\pi i+(-\infty,-R_\eps-1])=\emptyset$
is at least $1-\eps$.
If for any $x\le x_0$, $x+\pi i\in Q$, then $\psi(T,x+\pi i)-\xi(T)\in Q_T$.
Since $x\le x_0$, so $\Ree\psi(T,x+\pi i)-\xi(T)\le X(T)\le
-R_\eps-1$, and so $\psi(T,x+\pi i)-\xi(T)\in Q_T\cap(\pi i+(-\infty,-R_\eps-1])$.
Thus the probability that $Q\cap(\pi i +(-\infty,x_0])=\emptyset$
is at least $1-\eps$. Since $\eps>0$ is arbitrary, so a.s.\
$Q\cap(\pi i +(-\infty,x_0])=\emptyset$. Since this holds for any $x_0\in\N$, so a.s.\
$Q\cap\R_\pi=\emptyset$.

Fix $\eps>0$ and  $x_0\ge R_\eps+1$.
Let $X_0(t)=\psi(t,x_0)-\xi(t)$, $0\le t<T_0$, where $[0,T_0)$ is the largest interval on which
$\psi(t,x_0)$ is defined. Then $X_0(t)$ satisfies
the SDE: $dX_0(t)=-\sqrt\kappa dB(t)+\coth_2(X_0(t))dt-\sigma dt$.
Define $h_0$ on $(0,\infty)$ such that
$$h_0'(x)=\exp(x/2)^{\frac 4\kappa\,\sigma}\sinh_2(x)^{-\frac 4\kappa},\quad 0<x<\infty.$$
Since $\kappa>4$ and $\sigma\ge 1$, so $h_0$ maps $(0,\infty)$ onto $(L,\infty)$ for some $L\in\R$.
 From Ito's formula,  $Y_0(t):=h_0(X_0(t))$, $0\le t< T_0$, satisfies the SDE:
$dY_0(t)=-h_0'(X_0(t))\sqrt\kappa dB(t)$. Using a similar argument as before, we conclude that
a.s.\ $T_0<\infty$ and $\lim_{t\to T_0} X_0(t)=0$. So $T_0$ is a finite stopping time.
Let $\beta_{T_0}$ be the $\beta_T$ in Lemma \ref{betaT} with $T=T_0$.
Then $\beta_{T_0}$ has the same distribution as $\beta$.
Let $Q_{T_0}$ denote the set of subsequential limits
of $\beta_{T_0}(t)$ as $t\to\infty$. Then $Q_{T_0}=\psi(T_0,Q)-\xi(T_0)$.

Since $x_0$ is swallowed at time $T_0$, so $\xi(T_0)=d(T_0)$ and $b(T_0)\ge x_0$. Since the
extremal distance (c.f.\ \cite{Ahl}) between $(-\infty,a(T_0))$ and $(b(T_0),\infty)$ in $\St_\pi\sem L(T_0)$ is not
less than the extremal distance between them in $\St_\pi$, so from the properties of
$f_{T_0}$, we have $d(T_0)-c(T_0)\ge b(T_0)-a(T_0)$. Thus
$$c(T_0)-\xi(T_0)=c(T_0)-d(T_0)\le a(T_0) -b(T_0)\le -b(T_0)\le -x_0\le -R_\eps-1.$$
 If  $Q\cap(-\infty,a(T_0)]\ne\emptyset$,
then since $Q_{T_0}=\psi(T_0,Q)-\xi(T_0)$, so $Q_{T_0}\cap(-\infty,c(T_0)-\xi(T_0)]\ne\emptyset$, which
happens with probability less than $\eps$ since $\beta_{T_0}$ has the same distribution as $\beta$, and
 $c(T_0)-\xi(T_0)\le -R_\eps-1$. From Lemma \ref{boundaryofstrip}, for every $x\in(a(T_0),b(T_0))$,
 $L(T_0)$ contains a neighborhood of
$x$ in $\St_\pi$. Since $\beta$ does not cross its past, so $Q\cap(a(T_0),b(T_0))=\emptyset$.
Thus the probability that $Q\cap (-\infty,b(T_0))\ne\emptyset$
is less than $\eps$. Since $b(T_0)\ge x_0$, and $x_0\ge R_\eps+1$ is
arbitrary, so the probability that $Q\cap\R\ne\emptyset$ is
less than $\eps$. Since $\eps>0$ is arbitrary, so a.s.\
$Q\cap\R=\emptyset$. $\Box$

\begin{Corollary}
Let $\kappa>4$ and $\rho\ge \kappa/2-2$. Suppose $\gamma_*(t)$, $0\le t<\infty$, is
a chordal SLE$(\kappa;\rho)$ trace started from $(0;1)$. Then a.s.\ $\gamma_*$ has
no subsequential limit on $\R$. \label{Corollary1}
\end{Corollary}
{\bf Proof.} This follows from the above lemma and Proposition \ref{coordinate2}. $\Box$

\begin{Theorem}
Let $\kappa>4$ and $\rho\ge \kappa/2-2$. Suppose $\gamma(t)$, $0\le t<\infty$, is
a chordal SLE$(\kappa;\rho)$ trace started from $(0;0^+)$ or $(0;0^-)$. Then
a.s.\ $\lim_{t\to\infty}\gamma(t)=\infty$. \label{0;0+}
\end{Theorem}
{\bf Proof.} By symmetry, we only need to consider the case that the trace is started
from $(0,0^+)$. Let $Q$ be the set of subsequential limits of $\gamma$.
From Proposition \ref{coordinate}, for any $a>0$, $(a\gamma(t))$ has the same distribution
as $(\gamma(a^2 t))$. Thus $aQ$ has the same distribution as $Q$ for any $a>0$.
To prove that a.s.\ $Q=\{\infty\}$, we suffice to show that a.s.\ $0\not\in Q$.

Let $\zeta(t)$ and $\vphi(t,\cdot)$, $0\le t<\infty$, be the driving function and chordal
Loewner maps for $\gamma$. Let $X(0)=0$ and
$X(t)=\vphi(t,0^-)-\zeta(t)$ for $t>0$. Then $(X(t)/\sqrt\kappa)$ is a Bessel process
with dimension $\frac2\kappa(2+\rho)+1\ge 2$. So a.s.\ $\limsup_{t\to\infty} X(t)=\infty$. Let $T$ be the
first time that $X(t)=1$. Then $T$ is a finite stopping time.
Let $\gamma_*(t)=\vphi(T,\gamma(T+t))-\zeta(T)$, $t\ge 0$. Then $\gamma_*$ is a chordal
SLE$(\kappa;\rho)$ trace started from $(0;1)$. From the last corollary, $\gamma_*$ has
no subsequential limit on $\R$.
Let $g_{T}=\vphi(T,\cdot)^{-1}$. Then
$g_T$ extends continuously to $\lin{\HH}$, and $\gamma(T+t)=g_T(\gamma_*(t)+\zeta(T))$. From
the property of $\vphi(T,\cdot)$, we have $g_T(z)=z+o(1)$ as $z\to\infty$,
so $g_T^{-1}(0)-\zeta(T)\subset\R$ is bounded.
If $0\in Q$, then $\gamma_*$ has a subsequential limit on $g_T^{-1}(0)-\zeta(T)\subset\R$,
which a.s.\ does not happen. Thus a.s.\ $0\not\in Q$. $\Box$

\begin{Corollary}
Let $\gamma_*$ be as in Corollary \ref{Corollary1}. Then a.s.\ $\lim_{t\to\infty}\gamma_*(t)=\infty$.
\end{Corollary}
{\bf Proof.} Let $\gamma$ be a chordal SLE$(\kappa;\rho)$ trace started from $(0;0^+)$.
Let $\zeta(t)$ and $\vphi(t,\cdot)$, $0\le t<\infty$, be the driving function and chordal
Loewner maps for $\gamma$. Let $X(0)=0$ and $X(t)=\vphi(t,0^-)-\zeta(t)$ for $t>0$.
Let $T$ be the first time that $X(t)=1$. Then $T$ is a finite stopping time.
Let $\gamma_1(t)=\vphi(T,\gamma(T+t))-\zeta(T)$, $t\ge 0$. Then $\gamma_1$ has the same distribution
as $\gamma_*$. Since a.s. $\lim_{t\to\infty}\gamma(t)=\infty$, so a.s.
$\lim_{t\to\infty}\gamma_1(t)=\infty$. Since $\gamma_1$ has the same distribution
as $\gamma_*$, so a.s.\ $\lim_{t\to\infty}\gamma_*(t)=\infty$. $\Box$

\begin{Theorem}
Proposition \ref{greater2} also holds for $\kappa>4$. \label{+-infty>4}
\end{Theorem}
{\bf Proof.} This follows from the above corollary and Proposition \ref{coordinate2}. $\Box$

\vskip 3mm

Let $\kappa>4$, $p_0=x_0+\pi i\in\R_\pi$, $\rho_+,\rho_-,\rho_0\in\R$, and $\rho_++\rho_-+\rho_0=\kappa-6$.
Let $\beta(t)$, $0\le t<\infty$, be a strip SLE$(\kappa;\rho_+,\rho_-,\rho_0)$ trace
started from $(0;+\infty,-\infty,p_0)$. Let $\xi(t)$, $\psi(T,\cdot)$, and $L(t)$, $0\le t<\infty$,
be the corresponding driving function, strip Loewner maps and hulls. Then there is some Brownian
motion $B(t)$ such that $\xi(t)$ satisfies the SDE:
$$d\xi(t)=\sqrt\kappa dB(t)-\frac{\rho_+-\rho_-}2\,dt-\frac{\rho_0}2\,\coth_2(\psi(t,p_0)-\xi(t))dt.$$
Let \BGE X(t)=\Ree\psi(t,p_0)-\xi(t),\quad 0\le t<\infty.\label{defofX}\EDE Then $X(t)$ satisfies the SDE:
$$dX(t)=-\sqrt\kappa dB(t)+\frac{\rho_+-\rho_-}2\,dt+
\Big(\frac\kappa 2-2-\frac{\rho_++\rho_-}2\Big)\tanh_2(X(t))dt.$$
Define $h$ on $\R$ such that
$$h'(x)=\exp(x/2)^{-\frac 4\kappa\cdot\frac{\rho_+-\rho_-}2}\cosh_2(x)
^{-\frac 4\kappa\cdot(\frac\kappa 2-2-\frac{\rho_++\rho_-}2)},
\quad x\in\R.$$
Let $Y(t)=h(X(t))$, $0\le t<\infty$.
From Ito's formula, $Y(t)$ satisfies the SDE: $dY(t)=-h'(X(t))\sqrt\kappa dB(t)$.
For $0\le t<\infty$, let $u(t)=\int_0^t\kappa h'(X(s))^2 ds$. Then
$Y(u^{-1}(t))$, $0\le t<u(\infty):=\sup u([0,\infty))$, is a partial Brownian motion.
The behavior of $X(t)$ as $t\to\infty$ depends on the values of $\rho_+$ and $\rho_-$.
Now we suppose that
$\rho_+,\rho_-\ge \kappa/2-2$. Then $h$ maps $\R$ onto $\R$.
If $u(\infty)<\infty$, then a.s.\
$Y(u^{-1}(t))$ is bounded on $[0,u(\infty))$, so $X(t)$ is bounded on $[0,\infty)$. This then implies
that $u'(t)$ is uniformly bounded below by a positive constant, and so $u(\infty)=\infty$, which is
a contradiction. Thus a.s.\ $u(\infty)=\infty$, and so $\limsup_{t\to u(\infty)} Y(u^{-1}(t))=\infty$ and
$\liminf_{t\to u(\infty)} Y(u^{-1}(t))=-\infty$, which implies that
$\limsup_{t\to \infty} X(t)=\infty$ and $\liminf_{t\to \infty} X(t)=-\infty$.

\begin{Lemma}
Let $\beta$ be as above. If $\rho_+,\rho_-\ge \kappa/2-2$,
then a.s.\ $\beta$ has no subsequential limit on $\R\cup\{+\infty,-\infty\}
\cup\R_\pi\sem\{p_0\}$. \label{supinf}
\end{Lemma}
{\bf Proof.} Let $Q$ denote the set of subsequential limits of $\beta(t)$ as $t\to\infty$.
Let $L(\infty)=\cup_{t\ge 0} L(t)$.
From Theorem \ref{leftrightboundary} and Proposition \ref{coordinate2}, a.s.\ $p_0\in\lin{L(\infty)}$,
and $L(\infty)$
is bounded by two crosscuts in $\St_\pi$ that connect $p_0$ with a point on $(-\infty,0)$ and a
point on $(0,\infty)$, respectively. Thus a.s.\ $Q\cap(\R_\pi\cup\{+\infty,-\infty\}\sem\{p_0\})=\emptyset$.
Moreover, for any $\eps>0$, there is $R_\eps>0$ such that the probability that $\lin{L(\infty)}\cap\R\subset
[-R_\eps,R_\eps]$ is at least $1-\eps$.

For $r\in(0,1)$, let $A_r=\{z:r<|z-p_0|<\pi\}$. If $\dist(p_0,L(t))\le r$, then any curve
in $\St_\pi\sem L(t)$ that connects the arc $[p_0,+\infty)\subset\R_\pi$ with $(-\infty,a(t))$ must connect
the two boundary components of $A_r$. Thus the extremal distance between
$[p_0,+\infty)$ and $(-\infty,a(t))$ in $\St_\pi\sem L(t)$ is at least
$(\ln(\pi)-\ln(r))/\pi$. So the extremal distance between
$[\psi(t,p_0),+\infty)$ and $(-\infty,c(t))$ in $\St_\pi$ is at least
$(\ln(\pi)-\ln(r))/\pi$, which tends to $\infty$ as $r\to 0$. This implies
that $\Ree\psi(t,p_0)-c(t)\to \infty$ as $\dist(p_0,L(t))\to 0$. Similarly,
$d(t)-\Ree\psi(t,p_0)\to \infty$ as $\dist(p_0,L(t))\to 0$. Fix $\eps>0$. There is $r\in(0,1)$ such that if
$\dist(p_0,L(t))\le r$, then $\Ree\psi(t,p_0)-c(t),d(t)-\Ree\psi(t,p_0)\ge R_\eps+|x_0|+1$.
Let $T_0$ be the first $t$ such that $\dist(p_0,\beta(t))=r$. Since a.s.\ $p_0\in\lin{L(\infty)}$,
so $T_0$ is a finite stopping time.

Let $X(t)$ be defined as in (\ref{defofX}). Let $T$ be the first $t\ge T_0$ such that $X(t)=x_0=\Ree p_0$.
Since $\limsup_{t\to\infty}X(t)=+\infty$ and $\liminf_{t\to\infty} X(t)=-\infty$, so
 $T$ is also a finite stopping time.
Let $\beta_T$ be defined as in Lemma \ref{betaT}, then $\beta_T$ has the same distribution as $\beta$.
So the probability that $\lin{\beta_T([0,\infty))}\cap\R\subset [-R_\eps,R_\eps]$ is
at least $1-\eps$.
Since $\dist(p_0,L(T))\le \dist(p_0,L(T_0))=r$, so
$\Ree\psi(T,p_0)-c(T),d(T)-\Ree\psi(T,p_0)\ge R_\eps+|x_0|+1$. Since $X(T)=\Ree\psi(T,p_0)-\xi(T)
=x_0$, so $\xi(T)-c(T),d(T)-\xi(T)\ge R_\eps+1$, and so $[-R_\eps,R_\eps]\subset
[c(T)-\xi(T),d(T)-\xi(T)]$.
Thus the probability that $\lin{\beta_T([0,\infty))}\cap\R\subset [c(T)-\xi(T),d(T)-\xi(T)]$ is
at least $1-\eps$.
Since for every $x\in (a(T),b(T))$, $L(T)$ contains a neighborhood of $x$ in
$\St_\pi$, and $\beta$ does not cross its past, so $Q\cap(a(T),b(T))=\emptyset$.
If $Q\cap(-\infty,a(T)]\cup[b(T),\infty)\ne\emptyset$, then
$\beta_T$ has a subsequential limit on $(-\infty,c(T)-\xi(T)]\cup[d(T)-\xi(T),\infty)$,
which  happens with probability at most $\eps$. Thus the probability that $Q\cap\R\ne\emptyset$
is at most $\eps$. Since $\eps>0$ is arbitrary, so a.s.\ $Q\cap\R=\emptyset$.
$\Box$

\begin{Corollary}
Let $\kappa>4$, $\rho_+,\rho_-\ge \kappa/2-2$, and $p_-<0<p_+$. Let $\gamma_1(t)$, $0\le t<\infty$,
be a chordal SLE$(\kappa;\rho_+,\rho_-)$ trace started from $(0;p_+,p_-)$. Then a.s.\ $\gamma_1$ has
no subsequential limit on $\R$. \label{beta1}
\end{Corollary}
{\bf Proof.} This follows from the above lemma and Proposition \ref{coordinate2}. $\Box$

\begin{Theorem}
Let $\kappa>4$ and $\rho_+,\rho_-\ge \kappa/2-2$. Let $\gamma(t)$, $0\le t<\infty$,
be a chordal SLE$(\kappa;\rho_+,\rho_-)$ trace started from $(0;0^+,0^-)$. Then
a.s.\ $\lim_{t\to\infty}\gamma(t)=\infty$.\label{0+0-}
\end{Theorem}
{\bf Proof.} Let $Q$ be the set of subsequential limits of $\gamma$.
From Proposition \ref{coordinate}, for any $a>0$, $(a\gamma(t))$ has the same distribution
as $(\gamma(a^2 t))$. Thus $aQ$ has the same distribution as $Q$ for any $a>0$.
 So we suffice to show that a.s.\ $0\not\in Q$.

Let $\vphi(t,\cdot)$ and $\zeta(t)$ be the chordal
Loewner maps and driving function for the trace $\gamma$.
Then for $t>0$, $\vphi(t,0^-)<\zeta(t)<\vphi(t,0^+)$. Let $p_\pm =\vphi(1,0^\pm)-\zeta(1)$.
Let $\gamma_1(t)=\vphi(1,\gamma(1+t))-\zeta(1)$. Then conditioned on $\gamma(t)$, $0\le t\le 1$,
$\gamma_1$ is a chordal SLE$(\kappa;\rho_+,\rho_-)$ trace started from $(0;p_+,p_-)$.
From the argument in the proof of Theorem \ref{0;0+}, we see that if $0\in Q$, then
$\gamma_1$ has a subsequential limit on $\R$. From Corollary \ref{beta1}, this a.s.\
does not happen. Thus a.s.\ $0\not\in Q$. $\Box$

\begin{Theorem}
Let $\beta$ be as in Lemma \ref{supinf}. Then a.s.\ $\lim_{t\to\infty}\beta(t)=p_0$.\label{limitp0}
\end{Theorem}
{\bf Proof.} Let $\gamma(t)$, $0\le t<\infty$,
be a chordal SLE$(\kappa;\rho_+,\rho_-)$ trace started from $(0;0^+,0^-)$.
Let $\vphi(t,\cdot)$ and $\zeta(t)$ be the chordal
Loewner maps and driving function for the trace $\gamma$.
Let $\gamma_1(t)=\vphi(1,\gamma(1+t))-\zeta(1)$. Let $p_{\pm}=\vphi(1,0^\pm)-\zeta(1)$.
Then conditioned on $\gamma(t)$, $0\le t\le 1$,
$\gamma_1$ is a chordal SLE$(\kappa;\rho_+,\rho_-)$ trace started from $(0;p_+,p_{-})$. Choose
$W$ that maps $\HH$ conformally onto $\St_\pi$ such that $W(0)=0$ and $W(p_{\pm})=\pm\infty$.
Let $p_*=W(\infty)\in\R_\pi$, and $\rho_0=\kappa-6-\rho_+-\rho_-$.
From Proposition \ref{coordinate2}, there is a time-change function $u(t)$ such that
$\beta_*(t):=W(\gamma_1(u^{-1}(t)))$, $0\le t<\infty$, is a strip
SLE$(\kappa;\rho_+,\rho_-,\rho_0)$ trace started from $(0;+\infty,-\infty,p_*)$.
Let $\xi_*(t)$ and $\psi_*(t,\cdot)$, $0\le t<\infty$, denote the driving function and strip Loewner maps
for the trace $\beta_*$. Let $X_*(t)=\Ree\psi_*(t,p_*)-\xi_*(t)$, $0\le t<\infty$.
Let $T$ be the first time such that $X_*(t)=x_0=\Ree p_0$. Since
$\rho_+,\rho_-\ge\kappa/2-2$, so $\limsup_{t\to\infty}X_*(t)=\infty$ and
$\liminf_{t\to\infty}X_*(t)=-\infty$. Thus $T$ is a finite stopping time.
Let $\beta_T(t)=\psi_*(T,\beta_*(T+t))-\xi_*(T)$, $t\ge 0$. Then $\beta_T$ is a strip SLE$(\kappa;\rho_+,\rho_-)$
trace started from $(0;+\infty,-\infty,p_0)$. From Theorem \ref{0+0-}, we have a.s.\
 $\lim_{t\to\infty}\gamma(t)=\infty$,
which implies that $\lim_{t\to\infty}\gamma_1(t)=\infty$, and so $\lim_{t\to\infty}\beta_*(t)=p_*$.
Thus a.s.\ $\lim_{t\to\infty} \beta_T(t)=\psi(T,p_*)-\xi_*(T)=X_*(T)+\pi i=p_0$. Since
$(\beta_T(t))$ has the same distribution as $(\beta(t))$, so a.s.\ $\lim_{t\to\infty}\beta(t)=p_0$. $\Box$

\begin{Corollary}
Let $\gamma_1$ be as in Corollary \ref{beta1}. Then a.s.\ $lim_{t\to\infty}\gamma_1(t)=\infty$.
\end{Corollary}

\begin{Theorem} Proposition \ref{less2'} also holds for $\kappa>4$. \label{lim-1p}
\end{Theorem}
{\bf Proof}. This follows from Theorem \ref{Theorem-integration} and Theorem \ref{limitp0}. $\Box$

\begin{Theorem} Theorem \ref{lim-1p3} also holds for $\kappa>4$.
\end{Theorem}
{\bf Proof.} The proof of Theorem \ref{lim-1p3} still works here except that
Theorem \ref{lim-1p} should be used instead of
Proposition \ref{less2'}. $\Box$

\begin{Theorem}
Theorem \ref{tendstoinfty} also holds for $\kappa>4$. \label{tendstoinfty>4}
\end{Theorem}
{\bf Proof.} The proof of Theorem \ref{tendstoinfty} still works here except that
Theorem \ref{+-infty>4} and Theorem \ref{limitp0} should be used  instead of
Proposition \ref{greater2} and Proposition \ref{3pt**}. $\Box$

\vskip 3mm
Let $\gamma$ be as in Theorem \ref{tendstoinfty>4}. Let $K(t)$, $0\le t<\infty$, be the chordal Loewner
hulls generated by $\gamma$. Let $K(\infty)=\cup_{t\ge 0} K(t)$.
Let $\kappa'=16/\kappa$, $\rho_{\pm m}'=C_{\pm m}
(\kappa'-4)$, $1\le m\le N_\pm$, $\vec{\rho'}_\pm=(\rho'_{\pm 1},\dots,\rho'_{\pm N_\pm})$,
$C_\pm=\sum_{m=1}^{N_\pm}C_{\pm m}$, $W(z)=1/\lin{z}$, $p'_{\pm m}=W(p_{\pm m})$, $1\le m\le N_\pm$,
and $\vec{p'}_\pm=(p'_{\pm 1},\dots, p'_{\pm N_\pm})$.
In Lemma \ref{K-infty}, if we take $N_\mp+1$ force points, one of which is $x_1^+$, on $(x_1,x_2)$,
and take $N_\pm+1$ force points, one of which is $x_1^-$, outside $[x_1,x_2]$, then
we have the following theorem.

\begin{Theorem}
(i) If $N_+\ge 1$, then $W(\pa_\HH^+ K(\infty))$ has the same distribution
as a chordal SLE$(\kappa';(1-C_+)(\kappa'-4),(1/2-C_-)(\kappa'-4),\vec{\rho'}_+,
\vec{\rho'}_-)$ trace started from $(0;0^+,0^-,\vec{p'}_+,\vec{p'}_-)$. And
$\pa_\HH^+ K(\infty)$ is a crosscut in $\HH$ that connects $\infty$ with
some point that lies on $(0,p_1)$.
(ii) If $N_-\ge 1$, then $W(\pa_\HH^- K(\infty))$ has the same distribution
as a chordal SLE$(\kappa';(1/2-C_+)(\kappa'-4),(1-C_-)(\kappa'-4),\vec{\rho'}_+,
\vec{\rho'}_-)$ trace started from $(0;0^+,0^-,\vec{p'}_+,\vec{p'}_-)$. And
$\pa_\HH^- K(\infty)$ is a crosscut in $\HH$ that connects $\infty$ with
some point that lies on $(p_{-1},0)$.
\end{Theorem}

Let $\beta(t)$, $X(t)$, and $h(x)$  be defined as before Lemma \ref{supinf}.
Then $(h(X(t)))$ is a local martingale.
Let $I_1=[\kappa/2-2,\infty)$, $I_2=(\kappa/2-4,\kappa/2-2)$, and
$I_3=(-\infty,\kappa/2-4]$. Let Case (jk) denote the case that
 $\rho_+\in I_j$ and $\rho_-\in I_k$. We have studied Case (11).
In Cases (12) and (13), $h$ maps $\R$ onto $(-\infty,L)$ for some $L\in\R$, and we
conclude that a.s. $\lim_{t\to\infty}X(t)=\infty$. Symmetrically, in Cases
(21) and (31), a.s. $\lim_{t\to\infty}X(t)=\infty$. In Cases (22), (23), (32)
and (33), $h$ maps $\R$ onto $(L_1,L_2)$ for some $L_1<L_2\in\R$, and we
conclude that for some $p\in(0,1)$, with probability $p$, $\lim_{t\to\infty}X(t)=\infty$;
and with probability $1-p$, $\lim_{t\to\infty}X(t)=-\infty$.
Now we are able to prove the counterpart of Theorem 3.5 in \cite{duality}
when $\kappa>4$.

\begin{Theorem} In Case (11),
a.s.\ $\lim_{t\to\infty}\beta(t)=p_0$. In Case (12), a.s.\ $\lim_{t\to
\infty}\beta(t)\in (-\infty,p_0)$. In Case (21), a.s.\ $\lim_{t\to
\infty}\beta(t)\in (p_0,+\infty)$. In Case (13), a.s.\ $\lim_{t\to
\infty}\beta(t)=-\infty$. In Case (31), a.s.\ $\lim_{t\to
\infty}\beta(t)=+\infty$. In Case (22), a.s.\  $\lim_{t\to
\infty}\beta(t)\in (-\infty,p_0)$ or $\in (p_0,+\infty)$. In Case (23),
a.s.\ $\lim_{t\to \infty}\beta(t)=-\infty$ or $\in (p_0,+\infty)$. In
Case (32), a.s.\ $\lim_{t\to \infty}\beta(t)\in (-\infty,p_0)$ or
$=+\infty$. In Case (33), a.s.\ $\lim_{t\to \infty}\beta(t)=-\infty$
or $=+\infty$. And in each of the last four cases, both events
happen with some positive probability. \label{3pt***}
\end{Theorem}
{\bf Proof.} This follows from the same argument as in the proof of
Theorem 3.5 in \cite{duality} except that here we use
Theorem \ref{+-infty>4}, Theorem \ref{limitp0}, and Theorem \ref{lim-1p}. $\Box$

\vskip 3mm

We believe that for any chordal or strip SLE$(\kappa;\vec{\rho})$ trace $\beta(t)$,
$0\le t<T$, it is always true that a.s.\ $\lim_{t\to T}\beta(t)$ exists.


\begin{thebibliography}{99}
\bibitem{Ahl} Lars V.\ Ahlfors. {\it Conformal invariants: topics
in geometric function theory}. McGraw-Hill Book Co., New York, 1973.
\bibitem{dim-SLE} V.\ Beffara. The dimension of the SLE curves.
arXiv:math/0211322. To appear in {\it Ann.\ Probab.}.
\bibitem{cut-SLE} V.\ Beffara. Hausdorff dimensions for SLE$_6$,
{\it Ann.\ Probab.}, 32(3):2606-2629, 2004.
\bibitem{Julien-Duality} Julien Dub\'edat. Duality of Schramm-Loewner Evolutions,
arXiv:0711.1884.
\bibitem{LSW-2} Gregory F.\ Lawler, Oded Schramm and Wendelin Werner.
Conformal invariance of planar loop-erased random walks and uniform
spanning trees. {\it Ann. Probab.}, 32(1B):939-995, 2004.
\bibitem{RY} Daniel Revuz and Marc Yor. {\it Continuous Martingales
and Brownian Motion}. Springer-Verlag, 1991.
\bibitem{RS-basic} Steffen Rohde and Oded Schramm. Basic properties of
SLE. {\it Ann.\ Math.}, 161(2):883-924, 2005.
\bibitem{S-SLE} Oded Schramm. Scaling limits of loop-erased random walks
and uniform spanning trees. {\it Israel J. Math.}, 118:221-288,
2000.
\bibitem{duality} Dapeng Zhan. Duality of chordal SLE, arXiv:0712.0332v3. To appear
in {\it Inven.\ Math.}.
\bibitem{reversibility} Dapeng Zhan. Reversibility of chordal SLE,
arXiv:0705.1852. To appear in {\it Ann.\ Probab.}.
\bibitem{LERW} Dapeng Zhan. The Scaling Limits of
Planar LERW in Finitely Connected Domains. {\it Ann. Probab.}, 36(2):467-529, 2008.
\end{thebibliography}
\end{document}